\documentclass{article}
\usepackage{graphicx, subcaption, caption}
\usepackage{fullpage}

\usepackage{amsmath, amssymb, mathtools, bm}

\usepackage{booktabs, nicematrix}
\usepackage[export]{adjustbox}
\usepackage[colorlinks=true,linkcolor=blue]{hyperref}
\usepackage{cleveref}
\usepackage{xspace}
\usepackage[ruled,vlined]{algorithm2e}
\usepackage[numbers,sort]{natbib}
\usepackage{todonotes}

\DeclareMathOperator*{\argmin}{arg\,min}

\crefname{algocf}{alg.}{algs.}
\Crefname{algocf}{Algorithm}{Algorithms}

\crefformat{equation}{(#2#1#3)}
\crefrangeformat{equation}{(#3#1#4) to~(#5#2#6)}
\crefmultiformat{equation}{(#2#1#3)}{ and~(#2#1#3)}{, (#2#1#3)}{ and~(#2#1#3)}

\let\originalleft\left
\let\originalright\right
\renewcommand{\left}{\mathopen{}\mathclose\bgroup\originalleft}
\renewcommand{\right}{\aftergroup\egroup\originalright}

\newcommand{\review}[1]{#1}

\newcommand{\matern}{Mat\'ern\xspace}
\newcommand{\calC}{\mathcal{C}}
\newcommand{\calN}{\mathcal{N}}

\newcommand{\bbN}{\mathbb{N}}
\newcommand{\bbP}{\mathbb{P}}
\newcommand{\bbR}{\mathbb{R}}

\title{Deep Gaussian Process Priors for Bayesian Image Reconstruction}
\author{Jonas Latz\footnote{Department of Mathematics, University of Manchester, UK. \url{jonas.latz@manchester.ac.uk}} \and Aretha L.\ Teckentrup\footnote{School of Mathematics and Maxwell Institute of Mathematical Sciences, University of Edinburgh, UK. \url{a.teckentrup@ed.ac.uk}} \and Simon Urbainczyk\footnote{School of Mathematical and Computer Sciences and Maxwell Institute of Mathematical Sciences, Heriot-Watt University, Edinburgh, UK. \url{su2004@hw.ac.uk} (corresponding author)}}
\date{\today}

\begin{document}

\maketitle

\begin{abstract}
    In image reconstruction, an accurate quantification of uncertainty is of great importance for informed decision making. Here, the Bayesian approach to inverse problems can be used: the image is represented through a random function that incorporates prior information which is then updated through Bayes' formula. However, finding a prior is difficult, as images often exhibit non-stationary effects and multi-scale behaviour. Thus, usual Gaussian process priors are not suitable. Deep Gaussian processes, on the other hand, encode non-stationary behaviour in a natural way through their hierarchical structure. To apply Bayes' formula, one commonly employs a Markov chain Monte Carlo (MCMC) method. In the case of deep Gaussian processes, sampling is especially challenging in high dimensions: the associated covariance matrices are large, dense, and changing from sample to sample. A popular strategy towards decreasing computational complexity is to view Gaussian processes as the solutions to a fractional stochastic partial differential equation (SPDE). In this work, we investigate efficient computational strategies to solve the fractional SPDEs occurring in deep Gaussian process sampling, as well as MCMC algorithms to sample from the posterior. Namely, we combine rational approximation and a determinant-free sampling approach to achieve sampling via the fractional SPDE. We test our techniques in standard Bayesian image reconstruction problems: upsampling, edge detection, and computed tomography. In these examples, we show that choosing a non-stationary prior such as the deep GP over a stationary GP can improve the reconstruction. Moreover, our approach enables us to compare results for a range of fractional and non-fractional regularity parameter values.
\end{abstract}

\section{Introduction}
Gaussian processes (GPs) have become a popular tool to represent uncertainty in spatial objects, e.g., in hydrology \citep{sanchez2006representative,dodwell2015hierarchical}, climate modelling \citep{bolin2011spatial,lindgren11,sain2011spatial,sang2011covariance}, and medical imaging \citep{mejia2020bayesian,roininen2014whittle}.
GPs combine a number of desirable properties:
To start with, GPs are well understood theoretically. See, e.g., \citep{adler2010geometry,stein2012interpolation} as an introduction, as well as \citep{stuart2018posterior,teckentrup2020} for recent convergence results.
Moreover, they offer the flexibility to model a wide range of spatial phenomena. This flexibility is mainly due to the wide range of covariance functions that can be used to construct a Gaussian process.  
These covariance functions are often parametric and, thus, allow the user to set quantities such as correlation length, marginal variance, and regularity of the GP realisations directly.
Due to their theoretical simplicity and flexibility, they are especially popular as prior distributions in Bayesian inverse problems with linear forward operator and Gaussian noise.
In this specific case, the associated posterior distribution is given by a closed-form expression and does not rely on Monte Carlo sampling. The same closed-form expression is also used for GP regression. 

GPs present multiple computational challenges when employed in large-scale estimation problems. The closed-form expression used in GP regression requires the inversion of a matrix describing the observation covariance. This is often a computational bottleneck. In the context of Monte Carlo sampling, it is often necessary to decompose high-dimensional covariance matrices to sample approximations of GPs.
Different solution strategies have been proposed to mitigate these problems, depending on the application and the GP's covariance function. For instance, a common practice for representing \matern-type Gaussian processes is to see them as the solution to a stochastic partial differential equation (SPDE) \citep{whittle54,lindgren11}. This approach enables the use of state-of-the-art solvers for partial differential equations (PDEs) to speed up the evaluation of the GP regression mean and the sampling of the GP. The SPDE approach mainly profits from the locality of the differential operators and the sparsity of the associated system matrices. 
Further notable techniques to reduce computational effort include circulant embedding \citep{dietrich1997fast,istratuca2023smoothed}, adaptive cross-approximation \citep{harbrecht2012lowrank,harbrecht2015efficient,kressner20}, hierarchical matrices \citep{feischl2018fast,khoromskij09}, truncated Karhunen-Lo\`eve expansions \citep{schwab2006karhunen,saibaba2016randomized,contreras2018parallel}, sparse Gaussian processes \citep{titsias2009variational,smola2000sparse}, and Vecchia approximation \citep{katzfuss2021,sauer2023}.

The importance of uncertainty quantification has been recognised in a variety of different contexts, especially also in image reconstruction. Image reconstruction problems usually suffer from a combination of lack of data, observational noise, and instability. Bayesian inversion has recently gained popularity as an alternative to variational inversion that allows for accurate uncertainty quantification in the estimates.
Appropriate prior distributions in imaging are, e.g., total variation priors \citep{pragliola2023total} (although not discretisation-invariant \citep{lassas2004can}), Cauchy \citep{markkanen2019cauchy,suuronen2022cauchy} or other $\alpha$-stable priors \citep{suuronen2023bayesian,chada2019posterior}, \review{hierarchical Horseshoe and Student's $t$ priors \citep{dong2023inducing,uribe2023horseshoe,senchukova2024bayesian},} and machine-learning-based priors \citep{laumont2022bayesian,laumont2023maximum}. Image reconstruction problems are typically linear and, thus, Gaussian process priors are a natural choice \citep{roininen2014whittle,cohen1991,suuronen2022cauchy}. Unfortunately, commonly used GP models, such as \matern-type GPs, describe stationary behaviour in the object to be reconstructed and images of interest show multi-scale behaviour: they contain largely flat, constant areas as well as areas of small-scale details such as edges and texture, which, at the same time, are not captured well by a usual covariance function. 

\begin{figure}[htb]
    \centering
    \begin{subfigure}[t]{0.477\textwidth}
        \centering
        \includegraphics[height=0.5625\linewidth]{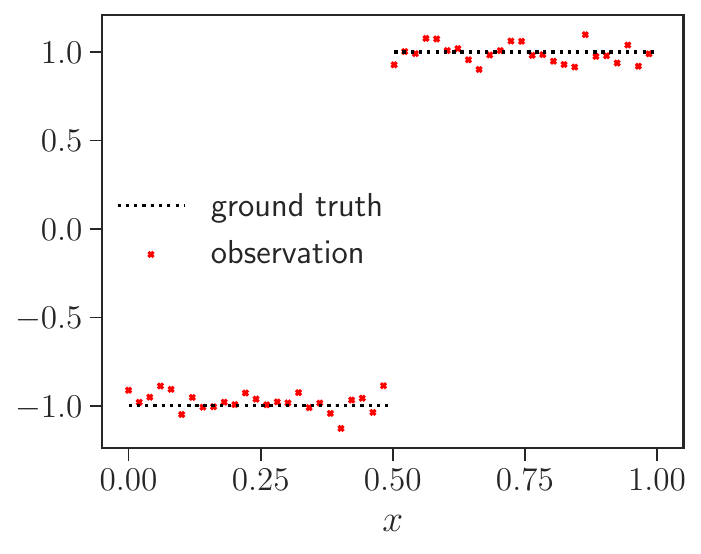}
        \caption*{Observations and ground truth.}
    \end{subfigure}
    \quad
    \begin{subfigure}[t]{0.477\textwidth}
        \centering
        \includegraphics[height=0.5625\linewidth]{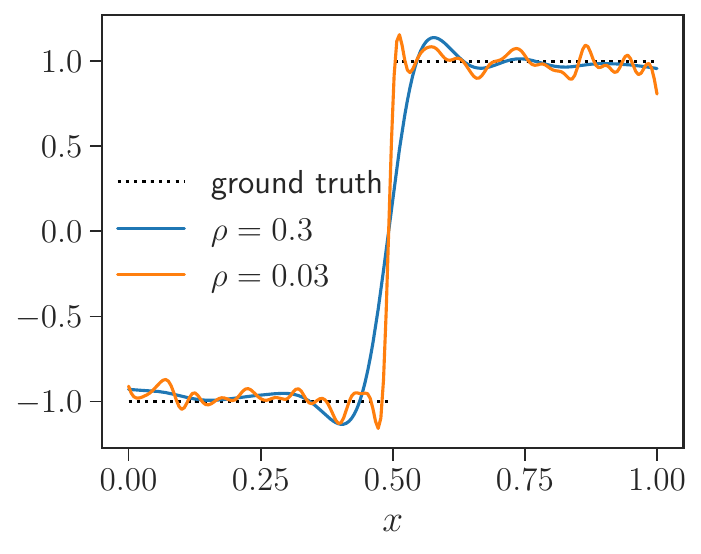}
        \caption*{GP regression mean for $\rho=0.03$ and $\rho=0.3$.}
    \end{subfigure}
    \caption{One-dimensional GP regression example with two correlation length values.}
    \label{fig:1d_example}
\end{figure}

We illustrate this last point in a small, one-dimensional example. The aim is to reconstruct the ground truth step function shown in the left plot in \Cref{fig:1d_example}. We have access to noisy observations of this ground truth, indicated as red crosses in the same plot. We perform Gaussian process regression using a \matern kernel \eqref{eq:matern} with $\nu=5/2$ and two different values for the correlation length: $\rho=0.03$ and $\rho=0.3$. The resulting regression mean functions are shown on the right. For $\rho=0.03$, we  capture the jump of the function at $x=0.5$ as a steep slope in our regression mean function. However, the noise that pollutes the observations is visible as oscillations in the actually flat area of the function. This makes sense intuitively, as the small value correlation length implies that we expect  small-scale features in the function. For the choice of the larger correlation length $\rho=0.3$, we make the reverse observation: The discontinuity at $x=0.5$ is very smooth, whereas  the regression mean is hardly affected by the observation noise reconstructing well the flat parts of the function. In this case, we would profit from a covariance function that can represent multiple length scales in different areas of the domain. These areas and the associated length scales are unknown and need to be identified as well.
In this work, we turn to \emph{deep Gaussian processes} \citep{dunlop18,damianou13,paciorek03,roininen2019hyperpriors} as a natural way to construct non-stationary spatial processes.


The denoising problem shown in Figure~\ref{fig:1d_example} is just one example of the image reconstruction problems we consider throughout this work. We more generally study image reconstruction problems that are given through associated inverse problems. We consider inverse problems of the form,
\begin{align}
    \label{eq:ip}
    A u^\dagger + \varepsilon^\dagger = d^\mathrm{obs},
\end{align}
where $u^\dagger$ is the \emph{unknown parameter or image} contained in some space $H \subseteq L^2(D)$, with some domain $D \subseteq \mathbb{R}^d$, $d^\mathrm{obs} \in \bbR^m$ is \emph{observed data}, $\varepsilon^\dagger\in \bbR^m$ is \emph{observational noise} and assumed to be a realisation of $\varepsilon \sim \mathcal{N}(0,\Gamma)$, with $\Gamma$ positive definite, and $A: H \rightarrow \mathbb{R}^m$ is a linear \emph{forward operator}.
Throughout this work, we consider the Bayesian approach to inverse problems. All random variables given throughout this work are defined on an underlying probability space $(\Omega, \mathcal{S}, \bbP)$. The Bayesian approach to inverse problems assumes that the unknown parameter is uncertain and represents the parameter by a random variable $U$ taking values in $H$. The random variable $U$ is independent of the noise $\varepsilon$ and it has a \emph{prior (distribution)} $\mathbb{P}(U \in \cdot)$ that represents the uncertainty in the parameter. In the Bayesian approach to inverse problems, we now incorporate the information in \eqref{eq:ip} into the distribution of $U$ conditioning on the data. We obtain the \emph{posterior (distribution)},
\begin{align*}
    \mathbb{P}(U \in \cdot \, | \, A U + \varepsilon = d^{\rm obs}) \: .
\end{align*}
Assuming that the prior distribution has a density $p(u)$ with respect to some reference measure, we obtain a density $p(\cdot|d^{\rm obs})$ of the posterior with respect to the same reference measure through Bayes' formula,
\begin{align}
    \label{eq:bayes}
    p\left(u \, | \, d^\mathrm{obs}\right)
    = \frac{p\left(d^\mathrm{obs} \, | \, u\right) \, p\left(u\right)}{p\left(d^\mathrm{obs}\right)} \: .
\end{align}
In the formulation of \cite{stuart2010inverse}, this reference measure is just the prior itself and $p(u) \equiv 1$. 
The \emph{likelihood} $p\left(d^\mathrm{obs} \, | \, u\right)$ is determined by the forward operator and the distribution of the noise. We have,
\begin{align}
    \label{eq:bayes_noise}
    p\left(d^\mathrm{obs} \, | \, u\right)
    &\propto \exp\left(- \frac{1}{2} \left\| \Gamma^{-1/2} \left(d^\mathrm{obs} - A u\right) \right\|^2\right) \: .
\end{align}
The \emph{model evidence} $p(d^{\rm obs})$ acts as a normalising constant and is usually difficult to evaluate. The posterior is usually approximated using Markov chain Monte Carlo (MCMC) methods that circumvent evaluations of $p(d^{\rm obs})$.

The focus of this work is the use of deep Gaussian process priors in Bayesian inverse problems of this form and associated computational strategies. Indeed, we propose a computational framework for deep Gaussian process priors in Bayesian image reconstruction and other linear inverse problems such as regression. In particular, our main contributions are the following:

\begin{itemize}
    \item We study spectral construction, finite element discretisation, and rational approximation of deep Gaussian processes constructed through the fractional SPDE formulation of \matern-type GPs -- extending previous works that only allowed for non-fractional SPDEs.
    \item Our computational framework utilises sparse computations within the SPDE framework. Notably, neither the forward operator $A$ nor the product $A^\ast A$ need to be sparse or stored in memory -- computations with them can be performed fully matrix-free. Overall, our methodology is scalable to high-dimensional images.
    \item We present MCMC methods to approximate posteriors with non-fractional and fractional Deep GP priors, the latter relying on a combination of marginalisation and normalisation-free methods for hyperpriors.
    \item We test and showcase our framework in extensive numerical experiments. These problems include image upsampling, edge detection, as well as computed tomography.
\end{itemize}

This work is organised as follows.
We introduce the class of deep Gaussian processes we consider throughout this work in \Cref{sec:deep_gp}. 
We explore the computational aspects of their SPDE representation in \Cref{sec:computations}. 
We describe and discuss the MCMC algorithms we use to solve inverse problems involving a deep Gaussian process prior in \Cref{sec:sampling}. 
We test our framework in numerical experiments in \Cref{sec:numerical}, and conclude our work in \Cref{sec:conclusion}.

\section{Deep Gaussian processes}
\label{sec:deep_gp}
Deep Gaussian processes are hierarchical extensions to standard Gaussian processes, which can model more complex behaviour.
There exists a number of different ways to construct a deep Gaussian process; for an overview we refer the reader to \citep{dunlop18}. In the rest of this section, we focus on a hierarchical construction where the length scale at each level (or ``layer'') is chosen dependent on the next-lower level. This enables us to model spatially varying length scales in a very intuitive manner.

First, however, we establish a bit of notation. In what follows, we indicate that a variable $U$ is a Gaussian process with mean function $m$ and covariance kernel $c$ by writing $U(x) \sim \mathrm{GP}(m(x), c(x, \cdot))$ for $x \in D \subseteq \bbR^d$.
In some places, we use the more general notation $U \sim \calN(m, \mathcal{C})$, where the covariance operator is defined via the integral identity,
\begin{align*}
    \mathcal{C} \varphi
    = \int_{\bbR^d} c(x, \cdot) \, \varphi(x) \: \mathrm{d}x 
    \quad
   \text{for } \varphi \in L^2\left(D\right) \: .
\end{align*}
We call this second notation more general, because it avoids the use of point values $U(x)$ and thus extends to the setting where the random variable $U$ takes values in a space of distributions.

\subsection{Gaussian processes and their SPDE representation}
\label{sec:gp_spde}
We base the construction of our deep Gaussian process on a Gaussian process with \matern covariance. This class of GPs is widely used in spatial statistics. They give the user access to individual parameters for marginal variance, length scale, as well as regularity of the GP realisations, allowing for great flexibility when modelling spatially correlated quantities.
A thorough overview of applications in which \matern covariances are used can be found in the introduction of \citep{khristenko19}.
The \matern covariance kernel 
is given by
\begin{align}
    \label{eq:matern}
    c(x,y) = \sigma^2 \mathcal{M}_v(\kappa \|x - y\|_2)
    \quad
    \text{for } x,y \in \bbR^d \: ,
    \quad
    \kappa = \frac{\sqrt{2 \nu}}{\rho} \: .
\end{align}
In the above, $\sigma ^ 2 > 0$ is the marginal variance, $\rho > 0$ denotes the length scale, and $\nu > 0$ is the smoothness, or regularity, of the resulting Gaussian process. Instead of the above, $\kappa$ is sometimes set to $\kappa = \sqrt{8\nu} / \rho$, see, e.g., \citep{lindgren11,croci18}. In both cases, $\kappa$ can be regarded as the (scaled) inverse length scale.
For ease of notation, we use the unit \matern function $\mathcal{M}_\nu$, which is defined as
\begin{align*}
\mathcal{M}_\nu(z) = \frac{z ^ \nu \: \mathcal{K}_\nu(z)}{2 ^ {\nu - 1} \: \Gamma(\nu)}
\quad
\text{for } z > 0 \: ,
\end{align*}
where $\mathcal{K}_\nu$ denotes the modified Bessel function of the second kind and order $\nu > 0$, and $\Gamma$ is the Gamma function. The exponential kernel $c(x,y) = \sigma^2 \exp(-\|x-y\|_2/\rho)$ is identical with the \matern kernel for $\nu = 1/2$, and the Gaussian kernel $c(x,y) = \sigma^2 \exp(-\|x-y\|_2^2/(2\rho^2))$ can be regarded as a limit of the \matern kernel for $\nu \rightarrow \infty$.

An advantage of using a \matern covariance is that the resulting Gaussian process can be seen as the solution to a stochastic partial differential equation \citep{whittle54,lindgren11}.
As we shall see later, this allows us to sample them efficiently.
Namely, \matern-type Gaussian processes are the stationary solutions to
\begin{align}
    \label{eq:spde}
    \left(\kappa^2 - \Delta\right) ^ {\alpha / 2} U(x, \omega) = \eta \, W(x, \omega)
    \quad
    \text{for } x \in \bbR^d \: , ~ \mathbb{P}\text{-a.e. } \omega \in \Omega \: .
\end{align}
Note that this SPDE is inherently defined on $\bbR^d$. Gaussian processes on a smaller domain can be obtained by restricting solutions to this SPDE to $D\subseteq\bbR^d$.
In the above, the regularity parameter $\alpha$ is determined by the value of $\nu$ as $\alpha = \nu + d/2$.
Furthermore, $\kappa$ is chosen as earlier in \cref{eq:matern}, and $\eta>0$ is a normalising constant,
\begin{align}
    \label{eq:normalisation}
    \eta ^ 2
    = \sigma ^ 2 \: \frac{\kappa ^ {2\nu} \: (4 \pi) ^ {d / 2} \: \Gamma\left(\nu + \frac{d}{2}\right)}{\Gamma(\nu)} \: .
\end{align}
This normalising constant is needed to ensure that the marginal variance of $U$ is equal to $\sigma ^ 2$.

The stochastic generalised function $W$ appearing on the right-hand side of \cref{eq:spde} is spatial white noise with unit variance. It can be seen as a generalised GP \citep{lindgren11} with the property that, for any set of sufficiently smooth test functions $\{\phi_1, \: \ldots, \: \phi_n\} \subset \mathcal{S} \subset L^2(\bbR^d)$, the pairings $(\langle\phi_i,W\rangle)_{1\leq i \leq n}$ are jointly normally distributed. For any $1\leq i,j \leq n$, the expectation and covariance of these random variables have to satisfy,
\begin{align*}
    \mathbb{E} \left[\langle\phi_i,W\rangle \right]
    &= 0 \: , \\
    \mathrm{Cov}\left(\langle\phi_i,W\rangle, \langle\phi_j,W\rangle\right)
    &= \langle\phi_i,\phi_j\rangle_{L^2(\bbR^d)} \: .
\end{align*}
Here, the pairing $\langle \cdot, \cdot \rangle$ denotes the evaluation of the linear functional $W$ w.r.t.\ a test function.
Since the right-hand side is a generalised function only, the stochastic partial differential equation (SPDE) in \cref{eq:spde} is to be understood in the sense of distributions only. The characterisation of the white noise used here naturally lends itself to be used in numerical methods based on weak formulations.

The use of the linear operator in \cref{eq:spde}, $(\kappa^2 - \Delta) ^ {\alpha / 2}$ requires some care. In the case that $\alpha / 2 \in \bbN$, this operator is well understood both in the analytical and computational sense. However, when $\alpha / 2 \not\in \bbN$, this operator is to be interpreted as a pseudo-differential operator \citep{krylov2008lectures,kwasnicki2017ten}.

We can now make use of \cref{eq:spde} to compute a sample of $U$ by generating a sample of the white noise, $W$, and then solving the SPDE for $U$.
For $\alpha/2\in\mathbb{N}$, discretisation of the pseudo-differential operator $(\kappa^2 - \Delta) ^ {\alpha / 2}$ typically results in a sparse matrix representation where the number of non-zero entries grows linearly with the number of spatial sample locations. Examples where such discretisations are used can be found in \citep{lindgren11,bolin2020rational,croci18}. Using an optimally preconditioned Krylov solver, the resulting linear system can be solved with a computational cost that is linear in the number of unknowns \citep{croci18,elman2014finite}.
This approach is very fast compared with most other sampling strategies for Gaussian processes. A na\"ive Cholesky decomposition of the covariance matrix, for instance, has a computational cost that is cubic in the number of unknown sample values.

In order to compute numerical solutions, one typically constrains the SPDE to a bounded connected domain $D \subset \bbR^d$. In the process, one has to introduce boundary conditions to ensure uniqueness of the solution. A common choice is to implement homogeneous Dirichlet or Neumann boundary conditions.
Choosing Neumann boundary conditions then additionally requires the boundary of our domain to be sufficiently smooth.
In practice, rather than \cref{eq:spde}, we thus consider the Whittle-\matern SPDE in a bounded domain,
\begin{equation}
    \label{eq:bounded_spde}
    \begin{aligned}
        \left(\kappa^2 - \Delta\right) ^ {\alpha / 2} U(x, \omega)
        = \eta \, W(x, \omega)
        \quad
        &\text{for } x \in D \: , ~ \mathbb{P}\text{-a.e. } \omega \in \Omega \: , \\
        \left.\begin{aligned}
            U(x, \omega) &= 0 \: , \\
            \text{or} \quad
            \frac{\partial U(x, \omega)}{\partial n} &= 0
        \end{aligned} \; \right\}
        \quad
        &\text{for } x \in \partial D \: , ~ \mathbb{P}\text{-a.e. } \omega \in \Omega \: .
    \end{aligned}
\end{equation}
Enforcing either Dirichlet or Neumann boundary conditions then leads to errors in the covariance of the SPDE solution as compared with \cref{eq:spde}. These errors will be more pronounced close to the boundary $\partial D$. A possible remedy is thus to increase the size of the computational domain. A thorough analysis of the error introduced by enforcing such artificial boundary conditions is described in \citep{khristenko19}.

\subsection{Deep Gaussian process construction}
Using the SPDE representation from \Cref{sec:gp_spde}, we now construct a deep Gaussian process.
To this end, consider a sequence of stochastic processes, $(U_\ell)_{\ell = 0, \ldots, L}$.
We say that this sequence represents a deep Gaussian process if, for each level $\ell=1, \ldots, L$, $U_\ell$ conditioned on $U_{\ell-1}$ is Gaussian and the covariance of $U_\ell$ depends on the next-lower level $U_{\ell-1}$.
Specifically, we choose the first process, or layer, $U_0$ as a Gaussian process. The process $U_0$ then determines the covariance kernel or operator of the next layer, $U_1$. This procedure is repeated until the top level $L$ is reached. The final layer $U_L$ then represents the output of the deep Gaussian process.

In order to model the dependence between layers, we use the SPDE approximation to \matern-type Gaussian processes in \cref{eq:bounded_spde}. While before we assumed the inverse length scale $\kappa$ to be a constant, we now choose $\kappa$ as a function of a realisation of the next-lower layer, $u_{\ell-1}(x)$. As a result, $\kappa$ is a spatially dependent function. Conditional on such a realisation of $U_{\ell-1}$, we then require that $U_\ell$ solves
\begin{align}
    \label{eq:deep_gp}
    \left(\kappa^2(u_{\ell-1}(x)) - \Delta\right) ^ {\alpha / 2} U_\ell(x, \omega) = \kappa^\nu(u_{\ell-1}(x)) \, \tilde{\eta} \, W_\ell(x, \omega)
    \quad
    \text{for } x \in D \: , ~ \mathbb{P}\text{-a.e. } \omega \in \Omega \: ,
    \quad
    + \text{ b.c.}
\end{align}
In the above equation, the scaling factor $\eta$ from \cref{eq:normalisation} is separated into $\kappa^\nu$ and a remaining factor of $\tilde{\eta} = \eta / \kappa^\nu$ to reflect the non-stationary nature of $\kappa$.

For a fixed regularity parameter $\alpha$ (or equivalently, $\nu$), the deep GP construction described here corresponds to controlling the (inverse) length scale through $u_{\ell-1}$.
This way, the resulting process $U_L$ can feature a multitude of locally varying length scales.
Except for the base level GP, $U_0$, all levels can model non-stationary behaviour without the need to explicitly specify a spatially dependent covariance function a priori.
This makes our deep Gaussian process an intuitively understandable and interpretable model for non-stationary quantities.

Alternatively, we can specify the distribution of $U_\ell$ in terms of a covariance operator $\mathcal{C}(u_{\ell-1})$ depending on a realisation of $U_{\ell-1}$.
To this end, we introduce $\mathcal{A}(u_{\ell-1})$ as a shorthand for the pseudo-differential operator appearing in \cref{eq:deep_gp},
\begin{align}
    \label{eq:operator}
    \mathcal{A}(u_{\ell-1})
    \coloneqq \tilde{\eta}^{-1} \kappa^{-\nu}(u_{\ell-1}) \, \left(\kappa^2(u_{\ell-1}(x)) - \Delta\right) ^ {\alpha / 2} \: .
\end{align}
Then, the construction in \cref{eq:deep_gp} tells us that $U_\ell$ conditioned on a realisation of $U_{\ell-1}$ follows a Gaussian distribution,
\begin{align}
    \label{eq:deep_normal}
    \bbP(U_\ell \in \cdot \: | \: U_{\ell-1} = u_{\ell-1})
    = \calN(0, \mathcal{C}(u_{\ell-1})) \: ,
    \quad
    \text{where } \mathcal{C}(u_{\ell-1}) ^ {-1} \coloneqq \mathcal{A}(u_{\ell-1}) ^ \ast \mathcal{A}(u_{\ell-1}) \: .
\end{align}
In the above, $\mathcal{A}(u_{\ell-1}) ^ \ast$ denotes the adjoint operator associated with $\mathcal{A}(u_{\ell-1})$.
The exact dependence of $\kappa$, and thus of the covariance operator $\mathcal{C}$, on $u_{\ell-1}$ is again a modelling choice. Following \citep{dunlop18}, we set
\begin{align}
    \label{eq:kappa}
    \kappa^2(u_{\ell-1}(x))
    = F(u_{\ell-1}(x)) \: ,
    \quad
    F(z) \coloneqq \min\left(F_- + a \cdot \exp(b \, z), ~ F_+\right) \: .
\end{align}
The reasoning behind choosing $F$ to be of this form is, on the one hand, that one can conveniently bound the values of $\kappa^2$ from above and below using the parameters $F_+ > F_- > 0$. Concretely, we need the lower bound to ensure that the linear operator defined this way is positive definite. On the other hand, the exponential growth term allows us to cover a range of length scales of multiple orders of magnitude, parameterised by $a$ and $b$. Any other choice of $F$ that is bounded in a similar way would be just as valid (compare with the assumptions on $\kappa$ in \citep[Assumption 1]{cox2020regularity}).

\begin{figure}[htb]
    \includegraphics[width=\linewidth]{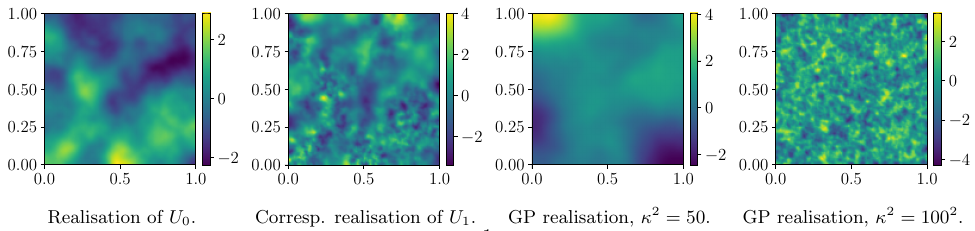}
    \caption{The left two plots show a realisation of a deep Gaussian process constructed using two layers, where the spatially varying length scale $\kappa^2 \in [50, 100^2]$. The right two plots show GP realisations corresponding to $\kappa^2=50$ and $\kappa^2=100^2$, respectively.}
    \label{fig:deep_gp}
\end{figure}

We present an illustration of this deep GP construction in \Cref{fig:deep_gp}.
For this example, we choose a grid of $128 \times 128$ sampling locations in $D = [0,1]^2$ and set $\alpha=3$. We select $\tilde{\eta}$ heuristically to achieve a marginal variance $\sigma^2 \approx 1$. The parameters of $F$ are chosen as $F_- = 50, ~ F_+ = 100^2, ~a = 1500, ~b = 1$.
A realisation of the bottom level GP, $u_0$, with constant $\kappa^2=200$ is shown on the far left, while we depict the corresponding realisation of the next-level GP next to it.
Areas in the left-most plot where $u_0$ assumes large values correspond to areas on the middle-left, where $u_1$ exhibits small length scales. Conversely, small values of $u_0$ mean small values of $\kappa(u_0)$ and thus a large local length scale of $u_1$. For comparison, we include plots of stationary \matern-type GP realisations with the inverse length scale given by $\kappa^2=F_-$, as well as $\kappa^2=F_+$. These two plots give a visual representation of the range of length scales that are covered by the deep GP construction.

In the above, one can interpret the parameters $F_-$ and $F_+$ as the lower and upper bounds for the inverse correlation lengths that can occur. The parameter $a$ sets the mode of $F(U_0)$, as the mode of the Gaussian process is the zero function. Finally, $b$ can be seen as a scaling parameter. In this two-layer example, changing $b$ is equivalent to directly scaling the marginal variance $\sigma^2$ of $U_0$.

\section{Computations with deep Gaussian processes}
\label{sec:computations}
As alluded to earlier, we aim at numerically solving the SPDE in \cref{eq:deep_gp}, satisfying some appropriate choice of boundary conditions.
Solving analytically, or simply evaluating a Karhunen--Lo\`eve expansion based on the eigenfunctions of the Laplacian as in \citep{khristenko19}, is not possible in our case since $\kappa$ is not constant. Instead, following \citep{lindgren11,bolin2020numerical}, we discretise the stochastic PDE using a finite element method. Using this discretisation then allows us to approximately compute realisations of a deep GP, as well as the covariance matrices for the Gaussian processes at each level.
In the remainder of this section, we first introduce a spectral definition of the operator in \cref{eq:deep_gp}. In the next step, we detail the finite element scheme used to discretise the operator, and finally give a rational approximation approach to deal with the fractional power in the SPDE representation of \matern-type GPs.

\subsection{Operator properties}
Before discretising the operator and function spaces relevant to the stochastic PDE in \cref{eq:deep_gp}, we take a step back and discuss their analytical properties in more detail. Specifically, we aim to clarify our use of the pseudo-differential operator.
Following \citep{bolin2020numerical,cox2020regularity}, we introduce a general linear operator, fulfilling some assumptions, as well as suitable function spaces. In a next step, we relate the functions and operators we use in the bounded Whittle-\matern SPDE \cref{eq:deep_gp} to these more general objects.

To this end, we denote with $H$ a separable Hilbert space and we let $\mathcal{L}: \mathcal{D}(\mathcal{L}) \subset H \rightarrow H$ be a self-adjoint, positive definite linear operator that is densely defined on $H$. Additionally, we assume that $\mathcal{L}$ has a compact inverse.
Then, $\mathcal{L}$ admits a system of eigenpairs $(\lambda_i, \phi_i)_{i\in\bbN}$, where the sequence $(\lambda_i)_{i\in\bbN}$ is non-decreasing, $\lambda_1 > 0$, and $\lim_{i \rightarrow \infty} \lambda_i = \infty$.
Furthermore, the eigenfunctions $\phi_i$ are orthonormal, that is, $\langle \phi_i, \phi_j \rangle_H = \delta_{i,j}$ for any $i,j\in\bbN$.
Using these eigenpairs, we can define the action of any positive power of $\mathcal{L}$ in the spectral sense \citep{lischke2020fractional},
\begin{align}
    \label{eq:spectral}
    \mathcal{L}^s u = \sum_{i\in\bbN} \lambda_i^s \langle \phi_i, u \rangle_H \, \phi_i
    \quad
    \text{for } s > 0 \: , ~ u \in \mathcal{D}(\mathcal{L}^s) \: .
\end{align}
For positive exponents $s>0$, the domain $\mathcal{D}(\mathcal{L}^s)$ of this fractional operator $\mathcal{L}^s$ is given as exactly the space where the above characterisation is well-defined,
\begin{align}
    \label{eq:dot_H}
    \dot{H}^{2s}
    \coloneqq \mathcal{D}(\mathcal{L}^s)
    = \left\{u \in H \, \Bigg| \, \sum_{i\in\bbN} \lambda_i^{2 s} \langle \phi_i, u \rangle_H^2 < \infty\right\}
    \quad
    \text{for } s > 0\: .
\end{align}
We identify $H = \dot{H}^0 \cong \dot{H}^{- 0}$ with its dual via its inner product using the Riesz representation theorem.
When $s < 0$, we denote with $\dot{H}^s = (\dot{H}^{- s})^\ast$ the dual w.r.t. $\langle\cdot,\cdot\rangle_H$ of $\dot{H}^{- s}$. The action of this negative-exponent operator is defined exactly as in \cref{eq:spectral}, with the exception that the inner product in $H$, $\langle\cdot,\cdot\rangle_H$, is replaced with the duality pairing $\langle\cdot,\cdot\rangle$ between $\dot{H}^{- s}$ and $\dot{H}^s$.
The operator $\mathcal{L}^s$ can then be seen as a one-to-one mapping between spaces $\dot{H}^t$ and $\dot{H}^{t - 2s}$ for any $s,t\in\bbR$. More precisely, there exists a unique continuous extension of $\mathcal{L}^s$ to an isometric isomorphism mapping $\dot{H}^t$ to $\dot{H}^{t - 2s}$ \citep[Lemma 2.1]{bolin2020numerical}.
As a result, for $f\in\dot{H}^{t - 2s}$, the equation $\mathcal{L}^s u = f$ admits a unique solution $u\in\dot{H}^t$.

\paragraph{Connection to Whittle--\matern SPDEs on bounded domains}
Now note that, with minimal assumptions on $\kappa$ and the domain $D$, we are in a setting where we can use the above analysis. Namely, we assume that $D\subset\bbR^d$ is bounded, connected and open.
Furthermore, we require that $\kappa\in L^\infty(D)$, which is automatically fulfilled due to our construction via $F$ in \cref{eq:kappa}.
For simplicity, we consider the case where we want to satisfy homogeneous Dirichlet boundary conditions in \cref{eq:deep_gp}. The Neumann case follows analogously, while requiring a stronger smoothness assumption on the boundary of $D$.
The Hilbert space we use in this setting is $H = L^2(D)$.

It now remains to specify the operator $\mathcal{L}=\kappa^2 - \Delta$ that is used in \cref{eq:deep_gp}.
The goal is to use this operator in the context of Sobolev spaces.
We make use of the Riesz representation to define $\mathcal{L}$ to uniquely characterise the action of $\mathcal{L}$ on any function $u \in H_0^1(D)$ via the identity,
\begin{align*}
    \langle\mathcal{L} u, v\rangle_H
    = \int_D \kappa^2 \, u \, v + \nabla u \cdot \nabla v \: \mathrm{d} x
    \quad
    \text{ for } v \in H_0^1(D) \: .
\end{align*}
The above integral is well-defined because $u,v\in H_0^1(D)$ and $\kappa\in L^\infty(D)$ is bounded.
Following the discussion in \citep{cox2020regularity}, the operator $\mathcal{L}$ introduced this way fulfils all the conditions to apply the analysis outlined earlier. Namely, $\mathcal{L}$ is defined on $H_0^1(D)$, which is dense in $H = L^2(D)$. Furthermore, $\mathcal{L}$ is self-adjoint and positive definite, while the inverse $\mathcal{L}^{-1}: H \rightarrow H$ is compact.

As a result, all the prerequisites to use the spectral operator construction \cref{eq:spectral} are fulfilled and we have a principled way of using $\mathcal{L}^s = (\kappa^2 - \Delta)^s$ in the context of the bounded Whittle-\matern SPDE.
In order to simplify notation, we write \cref{eq:deep_gp} as,
\begin{align}
    \label{eq:simple}
    \mathcal{L}^s U = B \: .
\end{align}
This is achieved by setting $B = \kappa^\nu \tilde{\eta} W$ and $s=\alpha / 2$, as well as using $\mathcal{L}$ as introduced in the above for either Dirichlet or Neumann conditions. For simplicity, we drop the subscripts relating to the deep GP level in \cref{eq:deep_gp}.
We now note that realisations $b$ of the right-hand side $B$ of \cref{eq:simple}, due to the spatial white noise term, almost surely have regularity $-d/2 - \varepsilon$ for any $\varepsilon>0$ (see, e.g., \citep[Proposition 2.3]{bolin2020numerical}).
That is, $b \in \dot{H}^{-d/2 - \varepsilon}$ and thus the corresponding realisation of the solution $u \in \dot{H}^{2s - d/2 - \varepsilon} \subset H = L^2(D)$ almost surely.
In fact, for our choice of operator $\mathcal{L}$ and $\mathcal{D}(\mathcal{L})$, the space $\dot{H}^{-d/2 - \varepsilon}$ coincides with the Sobolev space $H^{-d/2 - \varepsilon}$.
Since we can choose $\varepsilon>0$ such that $2s - d/2 - \varepsilon = \nu - \varepsilon > 0$, we have that $u\in H$ almost surely, and furthermore that $u$ has strictly positive regularity almost surely.

\subsection{Finite element discretisation of \matern SPDEs}
In order to compute numerical solutions to the fractional SPDE, \cref{eq:simple}, we introduce a family $(V_h)_{h > 0}$ of finite-dimensional subspaces of $\dot{H}^1$ (where $\dot{H}^1$ is given by \cref{eq:dot_H}). These subspaces could be defined, e.g., through a finite element discretisation. Denote with $n_h = \dim(V_h)$ the dimension of such a function space.
Using a Galerkin approximation approach, we introduce the discretised operator $\mathcal{L}_h: V_h \rightarrow V_h$ via the relation,
\begin{align*}
    \langle\mathcal{L}_h u_h , v_h\rangle_H = \langle\mathcal{L} u_h , v_h\rangle
    \quad
    \text{for } u_h,v_h\in V_h \: .
\end{align*}
We now turn to approximating the white noise term with an element in the finite-dimensional space, $V_h$.
Denote with $(e_{i;h})_{1 \leq i \leq n_h}$ an orthonormal basis of $V_h$. Then, we define the discretised white noise as,
\begin{align*}
    W_h = \sum_{i=1}^{n_h} \Xi_i \, e_{i;h} \: ,
    \quad
    \text{where } \Xi_i \sim \mathcal{N}(0,1) \text{ i.i.d.}
\end{align*}
It can easily be seen that in $V_h$, the discretised white noise is identical with white noise. That is, for any set of functions $(v_{1;h}, ~\ldots, ~v_{n;h}) \subset V_h$, the pairings $\langle W_h,v_{j;h}\rangle$ are jointly normally distributed and, for all $j,k\in\bbN$, we have
\begin{align*}
    \mathbb{E}\left[\langle W_h, v_{j;h}\rangle\right]
    &= \mathbb{E}\left[\sum_{i=1}^{n_h} \Xi_i \, \langle e_{i;h}, v_{j;h}\rangle\right]
    = 0 \: , \\
    \mathrm{Cov}\left(\langle W_h, v_{j;h}\rangle, \langle W_h, v_{k;h}\rangle\right)
    &= \mathbb{E}\left[\langle W_h, v_{j;h}\rangle \langle W_h, v_{k;h}\rangle\right]
    = \sum_{i=1}^{n_h} \mathbb{E}\left[\Xi_i^2\right] \langle e_{i;h}, v_{j;h}\rangle \langle e_{i;h}, v_{k;h}\rangle \\
    &= \langle v_{j;h}, v_{k;h}\rangle \: .
\end{align*}
While the second statement, concerning the covariance, is a necessary condition to ensure that we have constructed white noise, it also hints at how to obtain realisations of this discretised noise in practice.
In particular, it holds true when $(v_{1;h}, ~\ldots, ~v_{n_h;h}) \subset V_h$ are finite element basis functions.
In this setting, we can compute the mass matrix $M_h$, defined through $(M_h)_{j,k}=\langle v_{j;h}, v_{k;h}\rangle$. Let further $\bm w$ be the coefficient vector representing $W_h$ with respect to the basis $(v_{1;h}, ~\ldots, ~v_{n_h;h})$. Then, from the above statement we have that the covariance of the product $M_h\bm w$ is given by $M_h$, i.e., the covariance of $\bm w$ is the inverse mass matrix $M_h^{-1}$. A sample of the discretised noise can then be constructed as $\bm w = C^{-T} \bm\xi$, where $\bm\xi$ is drawn from the standard normal distribution $\calN(0, I_{n_h})$, and $C C^T = M_h$.

Coming back to the stochastic PDE, we denote with $\kappa_h^\nu$ an approximation to $\kappa^\nu$ in $V_h$. Then, the discretised equation that we want to solve reads,
\begin{align}
    \label{eq:discrete}
    \mathcal{L}_h^s U_h = \kappa_h^\nu \, \tilde{\eta} \, W_h \: .
\end{align}
In practice, one usually computes a realisation $w_h$ of $W_h$ and then evaluates the corresponding realisation of the solution, $u_h$.
The pointwise multiplication with $\kappa_h^\nu$ can simply be realised as multiplication with a diagonal matrix in a finite element setting. Solving with $\mathcal{L}_h^s$ turns out to be more challenging.

As $\mathcal{L}_h$ has finite rank, we can use an eigendecomposition of $\mathcal{L}_h$ to construct the inverse fractional operator $\mathcal{L}_h^{-s}$ (needed to solve \cref{eq:discrete}): When using the matrix representation of $\mathcal{L}_h$, denoted here by $L_h$, this corresponds to finding the eigendecomposition of said matrix.
Assume we have access to such a factorisation of the form $L_h = U \Lambda U^{-1}$. Then, the matrix representation of $\mathcal{L}_h^{-s}$  is given by $L_h^{-s} = U \Lambda^{-s} U^{-1}$. We will refer to this approach as the \textit{discrete eigenfunction method}, or DEM \citep{lischke2020fractional,hofreither2020unified}. The drawback of this method is the need to evaluate the eigendecomposition of $L_h$, which has a computational cost of order $\mathcal{O}(n_h^3)$. As a result, this approach is infeasible  for large $n_h$ and we propose to use a rational approximation instead.

\subsection{Rational approximation to fractional operator}
The idea behind rational approximation is to find a rational function $r(z)$ that is a good estimate of the (scalar) power function $z^{-s}$ \citep{hofreither2020unified,bolin2020rational,harizanov18}.
For now, assume that $0 < s < 1$.
The advantage of the rational function $r$ is that it can be applied to the operator $\mathcal{L}_h$ without the need of computing an eigendecomposition.
We usually expect that $r(z) \approx z^{-s}$ also implies that $r(\mathcal{L}_h)$ approximates $\mathcal{L}_h^{-s}$. 
In fact, when $L_h$ is diagonalisable, the rational function as used above can be seen as simply acting on the eigenvalues of $L_h$. To show this, let $U \Lambda U^{-1} = L_h$ denote the matrices resulting from diagonalising $L_h$. Then, as is the case for any matrix function, we have that,
\begin{align*}
    r(L_h) = r\left(U \Lambda U^{-1}\right) = U \, r(\Lambda) \, U^{-1} \: .
\end{align*}
As a result, when we say a ``good'' estimate, we mean that we want $r(z)$ and $z^{-s}$ to be uniformly close on an interval containing the eigenvalues of $\mathcal{L}_h$.
Note that the discrete operator $\mathcal{L}_h$ has a finite number of eigenvalues. Let us denote these by $0 < F_- \leq \lambda_{1;h}, \ldots, \lambda_{n_h;h} < \infty$. The uniform lower bound $F_-$ on the eigenvalues results from the construction of $\kappa$ in \cref{eq:kappa}.
We restrict the search for an optimal approximation $r$ to rational functions that are included in $\mathcal{R}_k$, the space of rational functions constructed from polynomials of maximum order $k\in\bbN$,
\begin{align*}
    \mathcal{R}_k = \left\{r = \frac{p}{q} \; \bigg| \; p, q \in \mathcal{P}_k \: , ~ q \not\equiv 0 \right\} \: .
\end{align*}
In summary, we choose our rational approximation function $r$ to be optimal in the following sense,
\begin{align}
    \label{eq:r_opt}
    r = \argmin_{\tilde{r} \in \mathcal{R}_k} \left\| \tilde{r}(z) - z ^ {-s} \right\|_{L^\infty\left(\left[F_-,\lambda_{n_h;h}\right]\right)} \: .
\end{align}
Technically, $r$ depends on $h$ through the largest eigenvalue $\lambda_{n_h;h}$ of the discrete operator. However, we omit this dependence for ease of notation.
In general, the exact function that fulfills the optimality condition in \cref{eq:r_opt} is not known analytically. However, there exist a few methods to numerically compute its coefficients \citep{hofreither2020unified}.
Of these methods, we have found the BRASIL algorithm proposed in \citep{hofreither2021rational} to be best suited for our purposes, as it optimises the condition \cref{eq:r_opt} directly.

Having determined the rational function itself, it is often useful to express it in partial fraction decomposition form.
In practice, the degrees of both its numerator and denominator are generally equal. As a result, we can write $r$ as,
\begin{align*}
    r(z) = c_0 + \sum_{j = 1}^k \frac{c_j}{z - d_j} \: .
\end{align*}
In practice, the partial fraction decomposition is preferred over other representations because of its numerical properties. It consists of comparatively few terms and, most importantly, only contains terms that are of order $-1$ in $z$.
In order to now approximate an operator instead of a scalar, we employ $r$ as a matrix function and apply it to the discrete operator, $L_h$. That is, we approximate both $L_h^{-s}$ and $u_h$ via,
\begin{align*}
    L_h^{-s}
    &\approx r(L_h)
    = c_0 \, I_{n_h} + \sum_{j = 1}^k c_j \left(L_h - d_j \, I_{n_h}\right) ^ {-1} \: , \\
    u_h
    &\approx u_{h,k}
    \coloneqq r(L_h) \: b_h
    = c_0 \, b_h + \sum_{j = 1}^k c_j \left(L_h - d_j \, I_{n_h}\right) ^ {-1} b_h \: ,
    \quad
    \text{where } b_h \coloneqq \kappa_h^\nu \, \tilde{\eta} \, w_h \: .
\end{align*}
In the above, we introduce $u_{h,k}$ to denote the solution obtained through rational approximation. The matrix $I_{n_h}$ is the identity matrix of dimension $n_h$.
Using the above formulas means that solving for the intractable fractional operator $L_h^s$ can be approximated by solving for matrices $L_h - d_j \, I_{n_h}$ a total of $k$ times instead. Depending on the discretisation used, the matrices $L_h - d_j \, I_{n_h}$ can usually be inverted using sparse matrix operations.
As a result, the approximate action of $L_h^s$ on the right-hand side $b_h$ can be computed with a cost of $\mathcal{O}(k n_h)$.

For cases where $s \geq 1$, it is common practice to split the operator $L_h^s$ into an ``integer power part'' and a ``fractional power part'', that is, $L_h^s = L_h^{\lfloor s \rfloor} L_h^{s - \lfloor s \rfloor}$ and $\lfloor s \rfloor = \max\{n \in \mathbb{Z}: n \leq s\}$. Then, one can solve for $L_h^{\lfloor s \rfloor}$ using standard finite element solvers and use a rational approximation to solve for $L_h^{s - \lfloor s \rfloor}$. This is particularly justified if, as in our case, the rational approximation algorithm relies on polynomials of equal degree in the numerator and denominator.

In practical computations, our discretisation of the SPDE \cref{eq:deep_gp} has the following form. We obtain the coefficient vector $\bm u$ representing $u_{h,k}$ in a finite element basis through the relation,
\begin{align}
    \label{eq:half_covariance}
    G_h \bm u &= \bm\xi \: , \\ \nonumber
    \text{where } G_h &= C^{T} D L_h^{\lfloor s \rfloor} r(L_h)^{-1} \quad \text{and } \bm\xi \text{ drawn from } \calN(0, I_{n_h}) \: .
\end{align}
In the above, C is chosen as a Cholesky factor of the mass matrix $M_h=CC^T$, and $D$ denotes a diagonal matrix representing multiplication with $\kappa_h^\nu$. The term $L_h^{\lfloor s \rfloor}$ is the integer power contribution of the operator $L_h^s$. When $s \notin \bbN$, we include the rational approximation term $r(L_h)^{-1}$ to represent $L_h^{s - \lfloor s \rfloor}$. In general, it is not practical to compute $G_h$ explicitly. In the case that $s=\alpha/2\in\bbN$, one can use mass lumping to obtain a sparse representation of $G_h$ \citep{lindgren11}. When $s\not\in\bbN$, we only have access to matrix-vector multiplications with $G_h$ or its inverse. The covariance matrix associated with $\bm u$ is $G_h^{-1}G_h^{-T}$. There are approaches to obtaining sparse representations of this covariance even when $s\not\in\bbN$ \citep{lindgren11,bolin2023covariance}, however, these cannot be applied in our case because of the presence of $D$, that is, the point-wise multiplication with $\kappa_h^\nu$.

\section{Sampling from posterior density}
\label{sec:sampling}
Having discussed the tools necessary for evaluating our deep GP prior, we now return to the inverse problem \cref{eq:ip} we aim to solve.
In the special case where the forward operator is linear, the noise is normally distributed, and the prior a Gaussian process, the posterior in \cref{eq:bayes} is known to be Gaussian as well \citep{kaipio2006statistical}.
In the discrete setting, there is an analytical formula for the mean vector and covariance matrix of such posterior distributions, see, e.g., \citep{williams2006gaussian}.
However, our goal is to use a deep Gaussian process prior. Since there is no closed-form expression of the posterior available in this case, we rely on Markov chain Monte Carlo (MCMC) sampling methods instead. MCMC methods construct an ergodic Markov chain that converges to the posterior distribution and, thus, allows us to estimate integrals with respect to the posterior measure.
Our MCMC algorithm of choice is a preconditioned Crank-Nicolson (pCN) method \citep{cotter2013mcmc,dunlop18}:
pCN algorithms are known to converge to the posterior distribution robustly with respect to the dimension of the discretised deep GP $n_h$ \citep{cotter2013mcmc,chen2018dimension}.
Note, however, that the robustness with respect to the dimension $n_h$ of the unknown parameter is proven only with the dimension $m$ of the observations fixed. Some dimension-dependence might arise in inverse problems where $m$ depends on $n_h$, as is often the case for imaging inverse problems.
In what follows, we describe the methods we use, as well as provide some computational strategies necessary for the implementation of such algorithms.

\subsection{Dimension-independent MCMC algorithms}
In order to be able to use a pCN-type proposal, and thus benefit of the dimension-robustness of the resulting sampler, we want to reformulate the posterior density in \cref{eq:bayes} to involve a Gaussian prior density.
One way to fulfil this prerequisite is to use ``whitened variables'' to realise the deep GP prior \citep{dunlop18,monterrubio2020posterior,chen2018dimension}. That is, we see the deep Gaussian process as a transformation of Gaussian white noise variables.
Such a transformation is readily at hand thanks to the SPDE representation: a \matern-type Gaussian process can be expressed as white noise that is transformed with a linear operator (in this case, a pseudo-differential operator). From the deep Gaussian process definition in \cref{eq:deep_gp,eq:operator}, we know that on each level $\ell = 1, \, \ldots , \, L$, and conditional on the realisation of the next-lower layer $u_{\ell-1}$, the layer $U_\ell$ of the deep GP is described by
\begin{align}
    \label{eq:whitened_layer}
    U_\ell
    = \mathcal{T}(u_{\ell-1}) W_\ell \: .
\end{align}
Here, we denote by $\mathcal{T}(u_{\ell-1})$ the solution operator depending on $u_{\ell-1}$, as detailed in \cref{eq:deep_gp}.
The white noise $W_\ell$ is indexed with $\ell$ to clarify that the white noise variables at each level are independent from each other.
Denote by $\overline{w} = (w_\ell)_{\ell = 0, \ldots, L}$ the vector containing realisations of the random quantities $W_\ell$, each of which follows a white noise distribution a priori.
We then introduce $\overline{\mathcal{T}}$ as the operator recovering the corresponding realisation of the deep GP layers, that is, the whole sequence $\overline{u} \coloneqq (u_\ell)_{\ell = 0, \ldots, L}$, from $\overline{w}$,
\begin{align}
    \label{eq:whitened}
    \overline{u} = \overline{\mathcal{T}}\left(\overline{w}\right) \: .
\end{align}
This relation means that, instead of sampling from the posterior of $(U_\ell)_{\ell = 0, \ldots, L}$, we can alternatively sample $(W_\ell)_{\ell = 0, \ldots, L}$.
The posterior density that we use to sample $(W_\ell)_{\ell = 0, \ldots, L}$ is then determined through Bayes' rule. That is,
\begin{align}
    \label{eq:bayes_whitened}
    p\left(\overline{w} \,|\, d^\mathrm{obs}\right)
    = \frac{p\left(d^\mathrm{obs} \,|\, \overline{w}\right) \: p\left(\overline{w}\right)}{p\left(d^\mathrm{obs}\right)} \: ,
\end{align}
where again $p\left(\overline{w}\right)$ and $p\left(\overline{w} \,|\, d^\mathrm{obs}\right)$ are densities with respect to an appropriate reference measure.
Effectively, we have rewritten the posterior density in \cref{eq:bayes} in terms of realisations of $(W_\ell)_{\ell = 0, \ldots, L}$, a quantity that is distributed with a Gaussian measure a priori.
Supposing that we have computed a chain $(\overline{w}_j)_{j\in\bbN}$ that samples $p\left(\overline{w} \,|\, d^\mathrm{obs}\right)$ to stationarity, we then obtain a chain $\left(\overline{\mathcal{T}}(\overline{w}_j)\right)_{j\in\bbN}$ that samples $p\left(\overline{u} \,|\, d^\mathrm{obs}\right)$ to stationarity simply by applying $\overline{\mathcal{T}}$ to each element of the original chain.

Through the forward problem, we have a direct correspondence between the observations $d^\mathrm{obs}$ and $\overline{u}$, which we are currently lacking for $\overline{w}$.
In order to make sense of the likelihood term in \cref{eq:bayes_whitened}, we thus use the knowledge that $p\left(d^\mathrm{obs} \,|\, \overline{w}\right) = p\left(d^\mathrm{obs} \,|\, \overline{u} = \overline{\mathcal{T}}\left(\overline{w}\right)\right)$.
Analogous to \cref{eq:bayes_noise}, we then obtain that,
\begin{align*}
    p\left(\overline{w} \,|\, d^\mathrm{obs}\right)
    \propto
    \exp\left(- \overline{\Phi}\left(\overline{w}; d^\mathrm{obs}\right)\right) \, p\left(\overline{w}\right) \: , \quad
    \text{where }
    \overline{\Phi}\left(\overline{w}; d^\mathrm{obs}\right)
    \coloneqq \frac{1}{2} \left\| \Gamma^{-1/2} \left(d^\mathrm{obs} - A\overline{\mathcal{T}}\left(\overline{w}\right)_L\right) \right\|^2 \: .
\end{align*}
In the above, we use $\overline{\mathcal{T}}\left(\overline{w}\right)_L$ to denote the top layer of the realisation $\overline{u}$, corresponding to $\overline{w}$.
Using the whitened representation, \cref{eq:whitened,eq:bayes_whitened}, comes with the advantage that we can now use a pCN proposal step \citep{cotter2013mcmc}, see also \Cref{alg:pcn}, when sampling from the posterior in \cref{eq:bayes_whitened}.


\subsection{Marginalising the top layer}
We now include an extra step to reduce the dimension of the state space that we sample our unknown in, as was done in \citep{dunlop18}.
In our experience, this significantly speeds up the convergence of the resulting Markov chain. We proceed as follows:

As the forward operator $A$ is linear, we have a closed-form expression \citep{williams2006gaussian} to describe the distribution of the top-most layer $U_L$, conditioned on $U_{L-1} = u_{L-1}$ and $A U_L + \varepsilon = d^\mathrm{obs}$.
That is, one can use standard formulas to obtain mean $m\left(u_{L-1};d^\mathrm{obs}\right)$ and covariance operator $\mathcal{C}\left(u_{L-1};d^\mathrm{obs}\right)$ such that,
\begin{align*}
    \bbP\left(U_L \in \cdot \: | \: U_{L-1} = u_{L-1}; d^\mathrm{obs}\right)
    = \calN\left(m\left(u_{L-1};d^\mathrm{obs}\right), \mathcal{C}\left(u_{L-1};d^\mathrm{obs}\right)\right) \: .
\end{align*}


As a result, it is only necessary to sample $(W_\ell)_{\ell = 0, \ldots, L}$ at layers up to index $L-1$.
These realisations of $(W_\ell)_{\ell\leq L-1}$, denoted with $\overline{w}_{<L}$, can then be transformed via \cref{eq:whitened_layer} to obtain realisations of $(U_\ell)_{\ell\leq L-1}$, while the top layer $U_L$ is represented by the GP regression mean computed using $u_{L-1}$.
We therefore update the posterior density to a marginal density referring to quantities on these sampled layers only,
\begin{align}
    \label{eq:bayes_marginalised}
    p\left(\overline{w}_{<L} \,|\, d^\mathrm{obs}\right)
    = \frac{p\left(d^\mathrm{obs} \,|\, \overline{w}_{<L}\right) \: p(\overline{w}_{<L})}{p(d^\mathrm{obs})} \: .
\end{align}
In order to make use of this identity in a sampling scheme, we need to be able to evaluate the likelihood, $p\left(d^\mathrm{obs} \,|\, \overline{w}_{<L}\right)$.
It is straightforward to derive the likelihood from the forward equation combined with the knowledge from \cref{eq:deep_normal} that $\bbP\left(U_L \in \cdot \: | \: U_{L-1} = u_{L-1}; d^\mathrm{obs}\right) = \mathcal{N}(0, \mathcal{C}(u_{L-1}))$.
Given a realisation $\overline{w}_{<L}$, and thus $u_{L-1}$, we can deduce that the observations in our forward problem \cref{eq:ip} follow a normal distribution with mean 0 and covariance $A \, \mathcal{C}(u_{L-1}) \, A^\ast + \Gamma$.
The quantity $d^\mathrm{obs}$, conditioned on $U_{L-1}=u_{L-1}$, is Gaussian and finite-dimensional, so that its conditional probability density exists. As a result, we obtain an expression for the likelihood $p\left(d^\mathrm{obs} \,|\, \overline{w}_{<L}\right)$ in \cref{eq:bayes_marginalised} that depends on $u_{L-1}$,
\begin{align}
    \label{eq:likelihood}
    \begin{aligned}
        p\left(d^\mathrm{obs} \,|\, \overline{w}_{<L}\right)
        &= p\left(d^\mathrm{obs} \,|\, u_{L-1}\right) \\
        &\propto \exp(- \Psi(u_{L-1}; d^\mathrm{obs})) \: , \\
        \text{where }
        \Psi(u_{L-1}; d^\mathrm{obs})
        &\coloneqq \frac{1}{2}\left(\|d^\mathrm{obs}\|^2_{\Sigma_{L-1}} + \log \det \Sigma_{L-1}\right) \: , \quad
        \Sigma_{L-1} \coloneqq A \, \mathcal{C}(u_{L-1}) \, A^\ast + \Gamma \: .
    \end{aligned}
\end{align}
The matrix norm in the above equation follows the convention $\|x\|_M^2 = x^T M^{-1} x$ for some positive definite matrix $M$ and vector $x$ of appropriate sizes.
Note that, even though $\mathcal{C}(u_{L-1})$ can be an infinite-dimensional operator in the above, $\Sigma_{L-1}$ is a positive definite matrix of size $m \times m$ by construction, where $m$ is the dimension of the observations, $d^\mathrm{obs}$.

\subsection{Sampling algorithm in the non-fractional case}
Using the results of the previous two subsections, we are now ready to state an algorithm that produces a Markov chain $(\overline{w}_j)_{j\in \bbN}$ with $p\left(\overline{w}_{<L} \,|\, d^\mathrm{obs}\right)$ as its stationary density.
Through the layer-wise transformation given in \cref{eq:whitened_layer}, we can then recover a chain $(\overline{u}_j)_{j\in\bbN}$ that samples the posterior of $(U_l)_{l\leq L-1}$ in stationarity.
The method that we use, \Cref{alg:pcn}, is based on the standard pCN method \citep{cotter2013mcmc}.
The likelihood used in \Cref{alg:pcn} is identical to \cref{eq:likelihood} and takes into account our earlier derivations for the whitened variables as well as the marginalised top layer.
The same algorithm is presented in \citep[Algorithm 1]{dunlop18}.

In the proposal step of the algorithm, we sample $\chi_j \sim \mathcal{N}(0, \mathcal{I})$. In the infinite-dimensional setting, the operator $\mathcal{I}$ denotes the identity operator and is chosen so that $\chi_j$ are samples of $(W_l)_{l\leq L-1}$. In finite dimensions, we select $\mathcal{I}$ to be the identity matrix of appropriate size. Then, $\chi_j$ are random vectors following a standard multivariate normal distribution.

\begin{algorithm}[ht]
    \DontPrintSemicolon
    Choose $\beta\in(0, 1]$, initial state $\overline{w}_0$ and number of steps $J\in\bbN$. \;
    
    \For{$j=0,~1,~\ldots,~J$}{
    Propose $\hat{w}_j = (1 - \beta^2)^{1/2}\overline{w}_j + \beta \chi_j$, where $\chi_j \sim \mathcal{N}(0, \mathcal{I})$. \;
    Set $\overline{w}_{j+1} = \hat{w}_j$ with probability $\alpha_j$,
    \begin{align*}
        \alpha_j
        &= \min \left\{1, \exp\left(\Psi\left(\upsilon_j; d^\mathrm{obs}\right) - \Psi\left(\hat{\upsilon}_j; d^\mathrm{obs}\right)\right)\right\} \: , \\
        \text{denoting }
        \upsilon_j
        &= \overline{\mathcal{T}}_{L-1}(\overline{w}_j)
        \text{ and }
        \hat{\upsilon}_j
        = \overline{\mathcal{T}}_{L-1}(\hat{w}_j)
        \: .
    \end{align*}
    
    Otherwise, set $\overline{w}_{j+1} = \overline{w}_j$. \;
    }
    \caption{MCMC algorithm with pCN proposals for sampling driving white noise.}
    \label{alg:pcn}
\end{algorithm}

In practice, \Cref{alg:pcn} will usually be used in a finite-dimensional setting. However, it remains valid in the infinite-dimensional case, where the random update steps $\chi_j$ are white noise samples.
It can sometimes be difficult to find a step size $\beta$ that allows for fast mixing of the chain. In order to avoid tuning it manually, it is common practice to choose it adaptively \citep{cotter2013mcmc}.
In our experiments, presented in \Cref{sec:numerical}, we choose the step size adaptively to achieve an acceptance rate of about 25\%.

\subsection{Likelihood computation in the non-fractional case}
\label{sec:computation_integer}
In order to make use of \Cref{alg:pcn}, we need to be able to compute the acceptance probabilities $\alpha_j$. These, in turn, rely on evaluations of the potential function $\Psi$, which is given in \cref{eq:likelihood}.
To obtain a first notion of how this potential can be computed, we restrict ourselves to the case where $\alpha / 2 \in \bbN$.
Then, after discretising, the inverse of the covariance operator $\calC(\upsilon_{L-1})$ can usually be written as a sparse precision matrix $P_h = G_h^T G_h$ \citep{lindgren11,dunlop18}.
As can be seen in \cref{eq:likelihood}, evaluations of $\Psi$ rely on computing the quadratic term $\|d^\mathrm{obs}\|^2_{\Sigma_{L-1}}$, as well as the logarithmic determinant $\log \det \Sigma_{L-1}$. In the discrete setting with a sparse precision matrix, the quadratic term is usually evaluated via the Woodbury matrix identity.
The sparsity of $P_h$ can then also be exploited to compute the logarithmic determinant by means of the identity $\log \det \Sigma_{L-1} = \log \det \left(P_h + A_h^\ast \Gamma^{-1} A_h\right) - \log \det \left(P_h\right) + \log \det \Gamma$, obtained via the cyclic property of the trace. Here, $A_h$ denotes a matrix representation of the forward operator $A$.
In our implementation, we obtain the quadratic term as well as the determinant through a Cholesky factorisations of $P_h + A_h^\ast \Gamma^{-1} A_h$ and $P_h$, and a subsequent solve. When the precision, as well as the product $A_h^\ast \Gamma^{-1} A_h$, can be represented with sparse matrices, this factorisation can be computed efficiently using fill-reducing permutations \citep{chen2008algorithm}. Note that, in order to use this approach, we not only require that $\alpha / 2 \in \bbN$, but also that the product of the forward operator with the noise covariance and itself, $A_h^\ast \Gamma^{-1} A_h$, is sparse.

\subsection{Sampling algorithm in the fractional case}
When describing \Cref{alg:pcn}, we observed that its efficient implementation hinges on the sparsity of the matrices $P_h$ and $A_h^\ast \Gamma^{-1} A_h$. However, when we do not have access to such sparse matrix representations, \Cref{alg:pcn} becomes infeasible already for moderate numbers of degrees of freedom.
This situation occurs when $\alpha/2 \not\in \bbN$, as we do not have access to an explicit precision matrix $P_h$ in this case.
Another scenario where using \Cref{alg:pcn} becomes impractical is when the forward operator leads to a matrix product $A_h^\ast \Gamma^{-1} A_h$ that is dense.
In these scenarios, we cannot make use of sparse direct methods to evaluate the potential function $\Psi$. While iterative solvers that rely on matrix-vector products only can still be used to compute the quadratic term in a reasonable time, in our experience the log-determinant remains a major issue. Even though a lot of progress has been made in recent years when it comes to trace estimation \citep{persson2022improved,cortinovis2022randomized,epperly2024xtrace}, a significant amount of matrix-vector products is needed to approximate the action of the logarithm of the operator, and then estimate its trace. As a result, generating a sufficiently precise estimate of the log-determinant makes the use of \Cref{alg:pcn} prohibitively expensive. Instead, we turn to determinant-free sampling methods. A number of algorithms have emerged that allow sampling in the case that the normalising determinant is intractable \citep{murray2012mcmc,ellam2017determinant}. Here, we make use of the methodology outlined in \citep{ellam2017determinant}. The advantage over the exchange algorithm proposed in \citep{murray2012mcmc} is that, in our setting, we require one iterative linear solve less.

The idea of the chosen determinant-free approach is to sample an auxiliary variable $Z$ alongside the desired samples of $(W_l)_{l\leq L-1}$, and ensure that the resulting marginal distribution of the samples of $(W_l)_{l\leq L-1}$ is identical to the target posterior. Specifically, we choose this auxiliary variable dependent on $U_{L-1}$,
\begin{align}
    \label{eq:aux_variable}
    \bbP(Z \in \cdot \: | \: U_{L-1} = u_{L-1})
    = \calN\left(0, \Sigma_{L-1}^{-1}\right) \: .
\end{align}
Since $U_{L-1}$ depends on $(W_l)_{l\leq L-1}$ deterministically, this gives rise to a posterior distribution for $Z$ when given observations $d^\mathrm{obs}$. We now aim to compute samples of $Z$ and $(W_l)_{l\leq L-1}$ from the joint posterior density,
\begin{align}
    \label{eq:joint_posterior}
    p\left(z, \, \overline{w}_{<L} \, | \, d^\mathrm{obs}\right)
    &\propto p\left(z \, | \, \overline{w}_{<L}\right) \, p\left(d^\mathrm{obs} \, | \, \overline{w}_{<L}\right) \, p\left(\overline{w}_{<L}\right) \: .
\end{align}
The conditional distribution of $Z$ was chosen such that the normalising determinants in the above likelihood term cancel out \citep{ellam2017determinant}. That is, the likelihood can be expressed without the need to evaluate $\det \Sigma_{L-1}$,
\begin{align}
    p\left(z \, | \, \overline{w}_{<L}\right) \, p\left(d^\mathrm{obs} \, | \, \overline{w}_{<L}\right)
    &= \mathrm{exp}\left(- \Phi\left(z, u_{L-1}; d^\mathrm{obs}\right)\right) \: , \nonumber \\
    \text{where } \Phi\left(z, u_{L-1}; d^\mathrm{obs}\right)
    &= \|z\|_{\Sigma_{L-1}^{-1}}^2 + \|d^\mathrm{obs}\|_{\Sigma_{L-1}}^2 \: .
    \label{eq:likelihood_rational}
\end{align}
In order to get MCMC samples $z_j$ and $\overline{w}_j$ from the joint posterior \cref{eq:joint_posterior}, we update both variables alternatingly in a Gibbs-like approach \citep{ellam2017determinant}. That is, for a given sample of $(W_l)_{l\leq L-1}$, we update the sample of $Z$ by generating an exact sample through the relation in \cref{eq:aux_variable}. With that new sample, we then take a single MCMC step based on the posterior \cref{eq:joint_posterior}. For this MCMC step, we again use a pCN proposal. The whole procedure is summarised in \Cref{alg:det_free}.
As earlier, we select the step size $\beta$ in \Cref{alg:det_free} adaptively to ensure an acceptance rate of about 25\%.

    
    

\begin{algorithm}[ht]
    \DontPrintSemicolon
    Choose $\beta\in(0, 1]$, initial state $\overline{w}_0$ and number of steps $J\in\bbN$. \;
    
    \For{$j=0,~1,~\ldots,~J$}{
    \tcp{Generate auxiliary variable $z_j$.}
    Compute a sample of $U_L$ via $v\sim\calN(0, \overline{\mathcal{T}}_{L-1}(\overline{w}_j))$. \;
    Sample observation noise, $e \sim \calN(0,\Gamma)$. \;
    Compute $z_j = \Sigma_{L-1}^{-1} (A v + e)$. \;
    \vspace{0.5\baselineskip}
    \tcp{Propose update $\hat{w}_j$.}
    Compute $\hat{w}_j = (1 - \beta^2)^{1/2}\overline{w}_j + \beta \chi_j$, where $\chi_j \sim \mathcal{N}(0, \mathcal{I})$. \;
    \vspace{0.5\baselineskip}
    \tcp{Accept/reject.}
    Set $\overline{w}_{j+1} = \hat{w}_j$ with probability $\alpha_j$,
    \begin{align*}
        \alpha_j
        &= \min \left\{1, ~ \exp\left(\Phi\left(z_j, \upsilon_j; d^\mathrm{obs}\right) - \Phi\left(z_j, \hat{\upsilon}_j; d^\mathrm{obs}\right)\right)\right\} \: , \\
        \text{denoting }
        \upsilon_j
        &= \overline{\mathcal{T}}_{L-1}(\overline{w}_j)
        \text{ and }
        \hat{\upsilon}_j
        = \overline{\mathcal{T}}_{L-1}(\hat{w}_j)
        \: .
    \end{align*}
    
    Otherwise, set $\overline{w}_{j+1} = \overline{w}_j$. \;
    }
    \caption{pCN-MCMC algorithm with auxiliary variable for sampling driving white noise.}
    \label{alg:det_free}
\end{algorithm}

We mentioned in passing that for each iteration of \Cref{alg:det_free}, we sample exactly from $\calN(0, \Sigma_{L-1}^{-1})$ to obtain our auxiliary variable. This is achieved by computing a sample of $U_L$, which we denote with $v$ for now.
This sample of the top layer is computed using the covariance operator induced by the newest available MCMC sample, $\overline{w}_j$.
Then, by simulating the forward problem, we have that $A v + e \sim \calN(0, \Sigma_{L-1})$. Here, $e$ denotes a sample of the observation noise in the forward problem \cref{eq:ip}. Finally, we solve with $\Sigma_{L-1}$ to obtain that $\Sigma_{L-1}^{-1} (A v + e) \sim \calN(0, \Sigma_{L-1}^{-1})$. That is, we need samples from $U_L$ and the observation noise, as well as a forward operator evaluation and a solve with $\Sigma_{L-1}$ to generate a sample of the auxiliary variable with the desired distribution. In practice, the solve with $\Sigma_{L-1}$ is the main bottleneck when sampling the auxiliary variable.

To be able to put \Cref{alg:det_free} to practical use, we need to specify how to evaluate the likelihood term appearing in the acceptance probability. This likelihood, as given in \cref{eq:likelihood_rational}, only contains quadratic terms. The norm $\|\cdot\|_{\Sigma_{L-1}^{-1}}$ only requires a multiplication with $\Sigma_{L-1}$ and is cheap to realise. The term containing $\|\cdot\|_{\Sigma_{L-1}}$, however, requires solving with $\Sigma_{L-1}$ and is computationally more involved. Since in the fractional setting we do not have access to the explicit entries of $\Sigma_{L-1}$, but can only evaluate matrix-vector products, it is natural to solve these systems iteratively. Our method of choice here is the LSQR solver \citep{paige1982lsqr,paige1982algorithm} applied to an equivalent least-squares problem.
\review{This allows us to solve with $\Sigma_{L-1}$ without forming this matrix explicitly, but instead only relying on $G_h$ as in \cref{eq:half_covariance}, the discretised forward operator $A_h$, and $\Gamma^{1/2}$.}


In order to reduce the number of iterations required to compute the solution to a certain accuracy, we employ a preconditioned variant of the LSQR algorithm, such as described in, e.g., \citep{arridge2014iterated}.
As a preconditioner, we propose to use $\Tilde{\Sigma}_{L-1}^{-1}$, where $\Tilde{\Sigma}_{L-1}$ denotes the covariance matrix associated with the next-largest value of $\alpha$ where this matrix has a sparse representation. That is, we choose the covariance matrix corresponding to the next-largest $\Tilde{\alpha}$ such that $\Tilde{\alpha}/2\in\mathbb{N}$. As in \Cref{sec:computation_integer}, we use a sparse Cholesky factorisation of this matrix in order to obtain solutions. In practice, one can usually leave $\Tilde{\Sigma}_{L-1}$ fixed for a number of steps in \Cref{alg:det_free}, as the MCMC updates to $\overline{w}_j$ tend to be small.

\section{Numerical results}
\label{sec:numerical}

In the preceding section, we have presented a framework that enables us to sample from the posterior distribution \cref{eq:bayes_marginalised} for any value of the parameter $\alpha > d/2$.
We now demonstrate the performance of our algorithm in a number of experiments.\footnote{\review{Code for all numerical experiments discussed here can be found at: \href{https://github.com/surbainczyk/deep_gp_priors}{https://github.com/surbainczyk/deep\_gp\_priors}}} These experiments are designed to highlight some possible use cases for deep GP regression. We also aim to give the reader an intuition for some practical considerations that arise when working with deep GPs in the proposed framework. In order to showcase the benefits of employing a non-stationary model for the prior, we compare the results of our framework with those obtained through standard GP regression.

\subsection{Image reconstruction with observation operator}
\label{sec:observation}

\begin{figure}
    \centering
    \begin{subfigure}[b]{0.25\textwidth}
        \centering
        \includegraphics[height=0.8\linewidth]{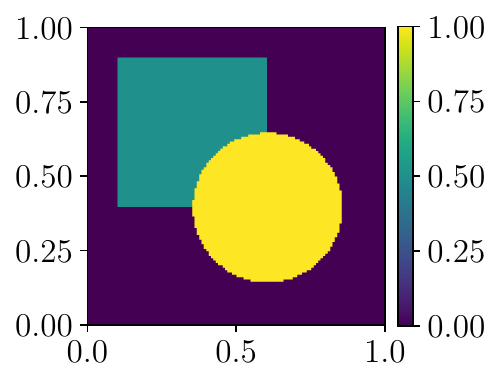}
        \caption*{``Square/circle'' image.}
    \end{subfigure}%
    \qquad
    \begin{subfigure}[b]{0.25\textwidth}
        \centering
        \includegraphics[height=0.8\linewidth]{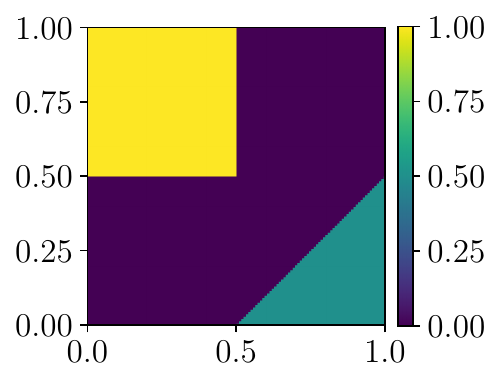}
        \caption*{``Corner/slope'' image.}
    \end{subfigure}%
    \qquad
    \begin{subfigure}[b]{0.25\textwidth}
        \centering
        \includegraphics[height=0.8\linewidth]{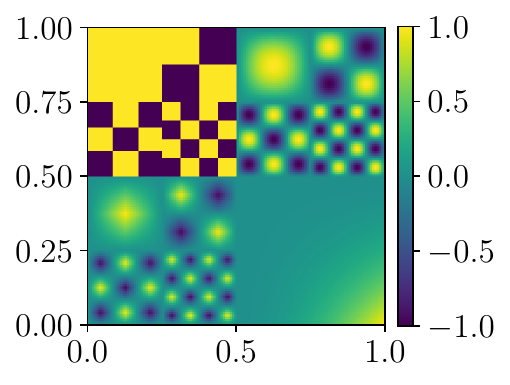}
        \caption*{``Multiscale'' image.}
    \end{subfigure}%
    \caption{Ground truth functions for reconstruction.}
    \label{fig:ground_truth}
\end{figure}

As a first experiment, we reconstruct a ground truth image from a limited number of observations. We restrict ourselves to single-channel images, in theory this experiment could also be carried out on every channel of a multi-channel image. Here, we consider the two left-most images shown in \Cref{fig:ground_truth}. We discuss results for the third ground truth image in \Cref{fig:ground_truth} in \Cref{sec:add_results}. 
Both images considered here are composed of step functions with corners and edges at various angles.
The edges in combination with the constant areas of the step functions can be regarded as multi-scale features: a GP with a large length scale would represent the constant areas reliably even in the presence of noise, while the edges and corners require a much smaller length scale.
We generated these images on a $128\times 128$ grid. The forward operator $A$ in this case is a simple observation operator: we observe $1/16$ of the original $128\times 128$ pixels. The observation locations are chosen to form a grid to achieve homogeneous coverage of the whole domain. That is, in this experiment we effectively perform an upsampling task. The observations are polluted by noise drawn from a normal distribution with mean $0$ and standard deviation $0.02$. This setup is very similar to an experiment in \citep{dunlop18}, where the authors reconstruct a much smoother ground truth composed of scaled sine functions.
One addition we make here is that we normalise the observations. That is, we subtract the mean of $d^\mathrm{obs}$ and then divide by the standard deviation of $d^\mathrm{obs}$ to obtain normalised observations.
In postprocessing, we then have to apply the inverse transformation on the reconstructed image. The normalisation step allows us to reuse the same setup of the prior for any scaling of the ground truth image values.

We choose the deep GP prior as a two-layer deep GP, i.e., with one hidden layer only. This choice is based on \citep{dunlop18}, where the authors observe that adding more layers is only beneficial when there is enough information available in the form of observations. In experiments not reported here, the additional layers led to more effort in the Monte Carlo sampling without showing a significant advantage over the two-layer setting.
We set the length scale of this bottom layer to $\rho\approx 0.063$. This is achieved by setting $\kappa^2=1500$ for $\alpha=4$ and scaling this value accordingly for other choices of $\alpha$ to leave the corresponding length scale invariant. In order to obtain the value of $\kappa^2$ that corresponds to the same correlation length $\rho$, but for a different value of $\alpha$, we can use the relation given in \cref{eq:matern} and multiply with $(2\alpha - 2)/6$.
The amplitude parameter $\sigma$ is chosen via \cref{eq:normalisation} to ensure a marginal variance of about 1 on both layers. The parameters of the function $F$ relating both layers, see \cref{eq:kappa}, are set to $F_- = 50, ~ F_+ = 100^2, ~a = 200, ~b = 1$.
For $\alpha=4$, this corresponds to a smallest possible correlation length of 0.024, and a largest possible correlation length of 0.346. The function $F$ is again scaled with a factor of $(2\alpha - 2)/6$ in order to keep these values invariant w.r.t.\ changes in $\alpha$.
We represent the Gaussian processes on each layer using continuous piece-wise linear finite elements, with the number of degrees of freedom matching the number of pixels in the ground truth images.

\review{In our experiments, we assume that the true standard deviation of the noise is known. Alternatively, noise levels can be estimated jointly with other stationary GP parameters using, e.g., maximum likelihood estimation \citep{williams2006gaussian}, or cross-validation methods \citep{sundararajan1999predictive}. In our experience, when the noise variance is under-estimated, our deep GP model tends to overfit the observational data.}

Our goal is to compare a range of different values for the regularity parameter $\alpha$. Note that we are here trying to evaluate which of the values for $\alpha$ are appropriate in such upscaling problems. In practice, $\alpha$ is a modelling choice that depends on the particular application \citep{sang2011covariance,minasny2005matern}.
In our experiments, we choose $\alpha$ from the following possible values,
\begin{align*}
    \alpha \in \{1.5, ~2, ~2.5, ~3, ~3.5, ~4\} \: .
\end{align*}
The number of steps that we simulate the Markov chain depends on the choice of $\alpha$, as the regularity parameter affects the algorithm we are using. In the non-fractional case, we run \Cref{alg:pcn} for $2.25 \times 10^5$ iterations for $\alpha=2$, and $10^5$ iterations for $\alpha=4$, respectively. In the fractional case, i.e., $\alpha/2\not\in\bbN$, we limit \Cref{alg:det_free} to $5.5\times10^4$ iterations for $\alpha=1.5$, and $2\times10^4$ iterations for the remaining values of $\alpha$. In both cases, we leave the preconditioner used in the LSQR solver fixed for 100 accepted MCMC steps.
The difference in numbers of MCMC steps is to allow for a fairer comparison, as each iteration in \Cref{alg:det_free} is more expensive compared to \Cref{alg:pcn}.
Additionally, the covariance matrix $\Sigma_{L-1}$ associated with $\alpha=2$ has less non-zero entries compared to the covariance matrix for $\alpha=4$. This directly affects the computation time of \Cref{alg:pcn} for the non-fractional case, but also influences \Cref{alg:det_free} through our choice of preconditioner.
As a result of limiting the total iterations to these values, the computational budget was roughly identical in all MCMC runs.
For all values of $\alpha$, we discard the first $10^4$ iterations as burn-in.
As alluded to earlier, and in order to achieve comparable values of the correlation length for different values of $\alpha$, we rescale the function $F$ used to relate upper layer $U_1$ and bottom layer $U_0$.
The parameter $k$ controlling the accuracy of the rational approximation in the fractional case is set to $k=3$.

\begin{figure}[p]
    \centering
    \includegraphics[width=\textwidth]{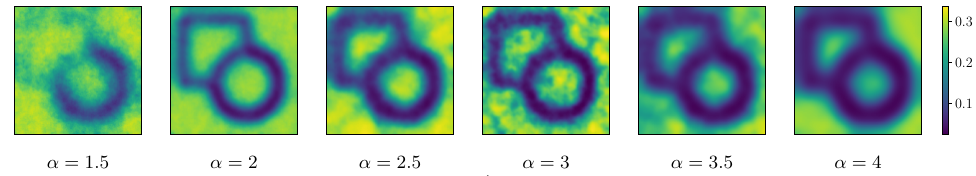}
    \caption{Correlation lengths used in deep GP reconstructions of ``square/circle'' image with varying $\alpha$.}
    \label{fig:corr_lengths_circle}
\end{figure}

\begin{figure}[p]
    \centering
    \includegraphics[width=\textwidth]{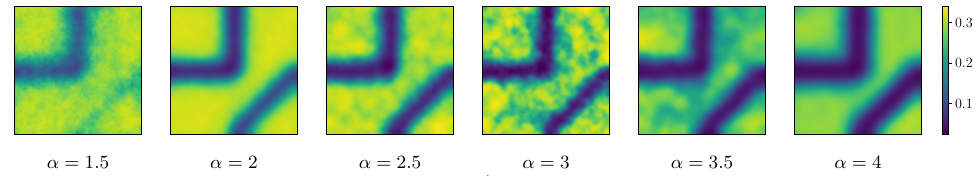}
    \caption{Correlation lengths used in deep GP reconstructions of ``corner/slope'' image with varying $\alpha$.}
    \label{fig:corr_lengths_corner}
\end{figure}

\paragraph{Discussion of results}
We begin exploring the reconstruction results by verifying how well the inner layer captures the structure of the image.
We give an overview of the estimated length scales for each value of $\alpha$, and both ground truth images, in \Cref{fig:corr_lengths_circle}, and \Cref{fig:corr_lengths_corner}, respectively.
These length scales were obtained in the following manner. Using \Cref{alg:pcn,alg:det_free}, we compute samples of the hidden layer $U_0$. With these samples, we estimate the mean of $F(U_0)^{1/2}$. This corresponds to the value of $\kappa$ which in turn defines the correlation length of the top layer $U_1$.
Using the local correlation length values, we can then compute the reconstructed image for each $\alpha$ by evaluating the GP regression mean of $U_1$.

We can observe that for all values of $\alpha$, we reliably find the edges present in the image as increased values of the transformed hidden layer, which directly translates into small correlation length values. That is, in areas where edges are present, the hidden layer ensures that the local correlation length used in the top layer is small.
Conversely, in areas where the ground truth is constant, we observe large local correlation length values.
This is very much desired, as the local correlation length influences the visual smoothness of our regression result. In areas of large correlation length (and constant ground truth), the noise cannot be recovered by the reconstruction and is ``smoothed out''. Where the ground truth exhibits small-scale features, i.e., edges, we want the local correlation length to be small to be able to capture those edges.

The estimated correlation length plots in the fractional case, $\alpha=1.5, ~ 2.5, ~ 3, ~ 3.5$, look noisy compared to the non-fractional case. This can be explained with the larger maximum number of MCMC iterations we choose when setting $\alpha=2$ or $\alpha=4$.

Having investigated the behaviour of the inner layer, we now turn towards quantifying the reconstruction error of the top layer regression results.
For both ground truth images, we compute the $L^1$- and $L^2$-errors, as well as the structural similarity (SSIM) \citep{wang2004image} and peak signal-to noise ratio (PSNR) \citep{gonzalez2009digital}.
The $L^1$- and $L^2$-errors are optimal at 0, SSIM is optimal at 1, and PSNR is optimal at $+\infty$.
For different values of $\alpha$, the obtained errors are very similar. However, we will compare the performance for different values of $\alpha$ in more detail later in this section.

\begin{figure}[p]
    \centering
    \begin{subfigure}[b]{0.477\textwidth}
        \centering
        \includegraphics[height=0.5625\linewidth]{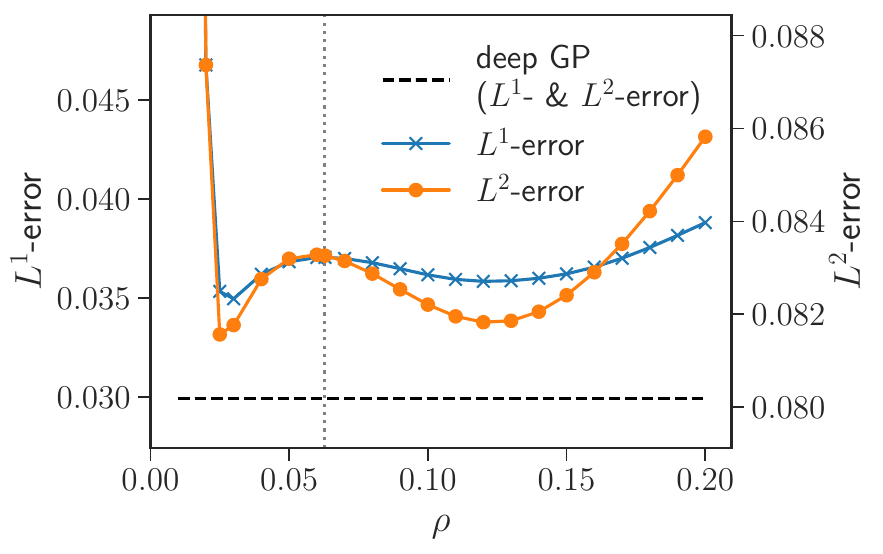}
        \caption*{$L^1$-errors and $L^2$-errors.}
    \end{subfigure}%
    \begin{subfigure}[b]{0.477\textwidth}
        \centering
        \includegraphics[height=0.5625\linewidth]{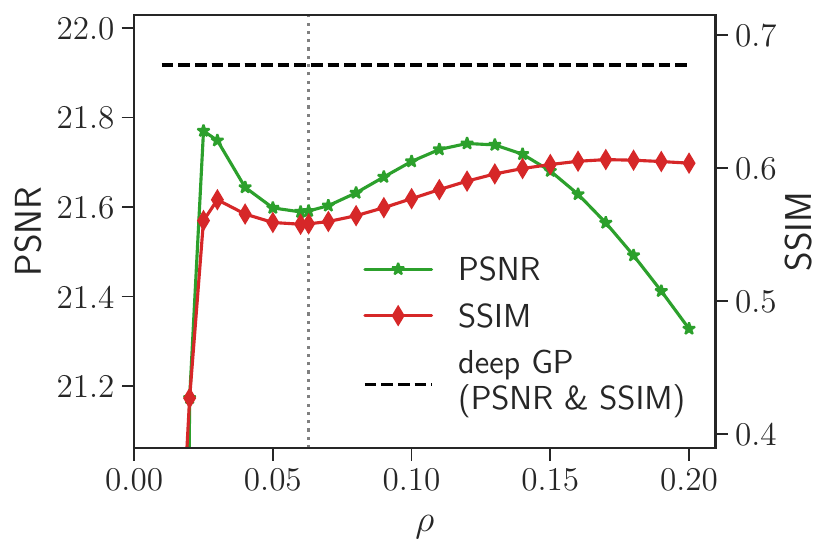}
        \caption*{PSNR and SSIM.}
    \end{subfigure}
    \caption{Reconstruction errors for ``square/circle'' image and $\alpha=3$, using standard GP and deep GP priors.}
    \label{fig:square_circle_errors_a3}
\end{figure}

\begin{figure}[p]
    \centering
    \begin{subfigure}[b]{0.477\textwidth}
        \centering
        \includegraphics[height=0.5625\linewidth]{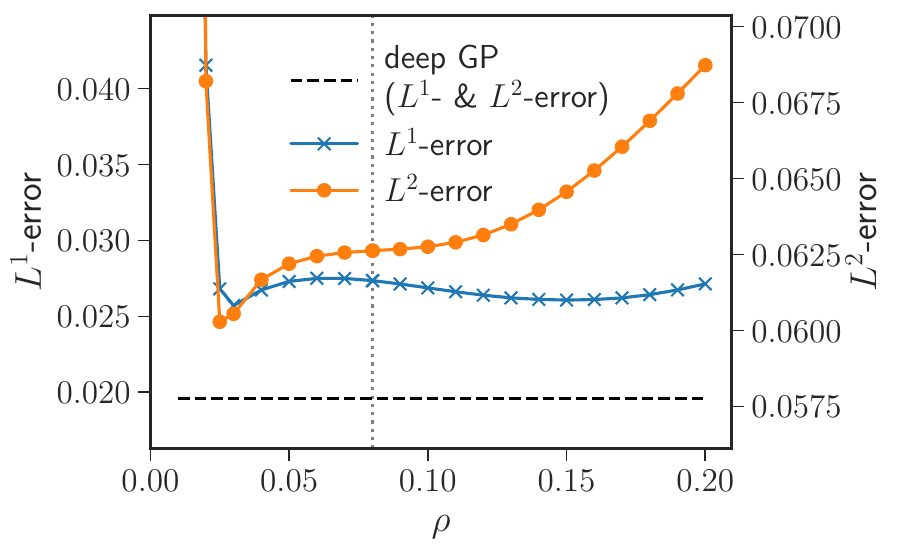}
        \caption*{$L^1$-errors and $L^2$-errors.}
    \end{subfigure}%
    \begin{subfigure}[b]{0.477\textwidth}
        \centering
        \includegraphics[height=0.5625\linewidth]{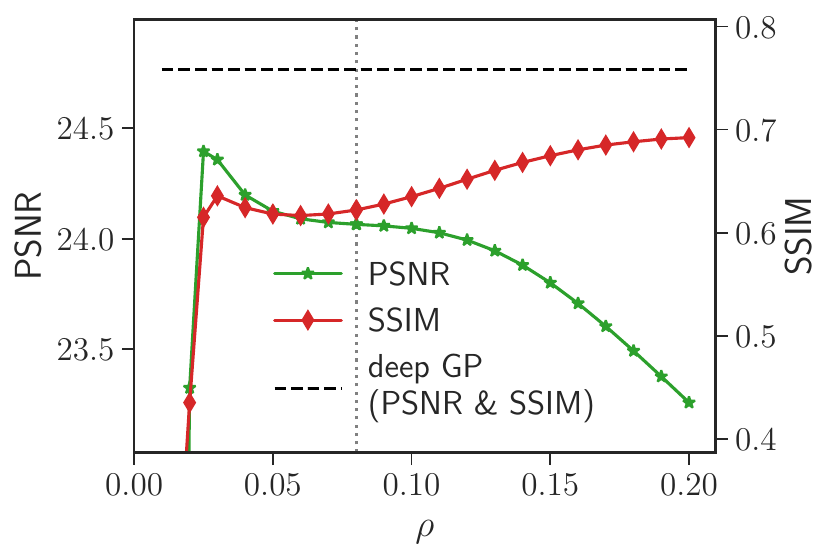}
        \caption*{PSNR and SSIM.}
    \end{subfigure}
    \caption{Reconstruction errors for ``corner/slope'' image and $\alpha=3$, using standard GP and deep GP priors.}
    \label{fig:corner_slope_errors_a3}
\end{figure}

On the other hand, we observe a larger benefit from using the deep GP setup, as compared to standard, stationary GP regression. As an example, we compare the errors obtained for a fixed regularity parameter, $\alpha=3$, and standard GP regression with a range of values for the correlation length $\rho$ with the flexible length scale deep GP prior. Additionally, we estimate an optimal length scale $\rho^\ast$ for GP regression by minimising the potential $\Psi$ from \cref{eq:likelihood} for a stationary GP covariance operator.
Evaluations of this potential are feasible for $\alpha=3$ since, in the stationary setting, we can use the approach from \citep{lindgren11} to compute a sparse precision matrix for the corresponding Gaussian processes.
The overviews for both ground truth images are presented in \Cref{fig:square_circle_errors_a3}, and \Cref{fig:corner_slope_errors_a3}, respectively. Lines for the $L^1$-error and $L^2$-error are given in the left plot of each figure, and the results for PSNR and SSIM are shown on the right. The performance of the deep GP is indicated by the horizontal dashed line. The estimated optimal length scale $\rho^\ast$ is given by the vertical dotted line in both plots. For both ground truth images, every choice of correlation length $\rho$, and each of the four error metrics, the deep GP outperforms the standard GP by a clear margin.

However, in order to observe a benefit from using the deep GP as a prior, it is crucial that the observations contain sufficient information to obtain a correlation length estimate that actually improves the reconstruction. In \Cref{sec:add_results}, we provide an example of a more complex ground truth image, where the deep Gaussian process fails to improve over the results obtained with stationary GP regression for some optimally chosen correlation length.

\begin{figure}
    \centering
    \begin{subfigure}[b]{0.5\textwidth}
        \raggedleft
        \includegraphics[height=0.536625\linewidth]{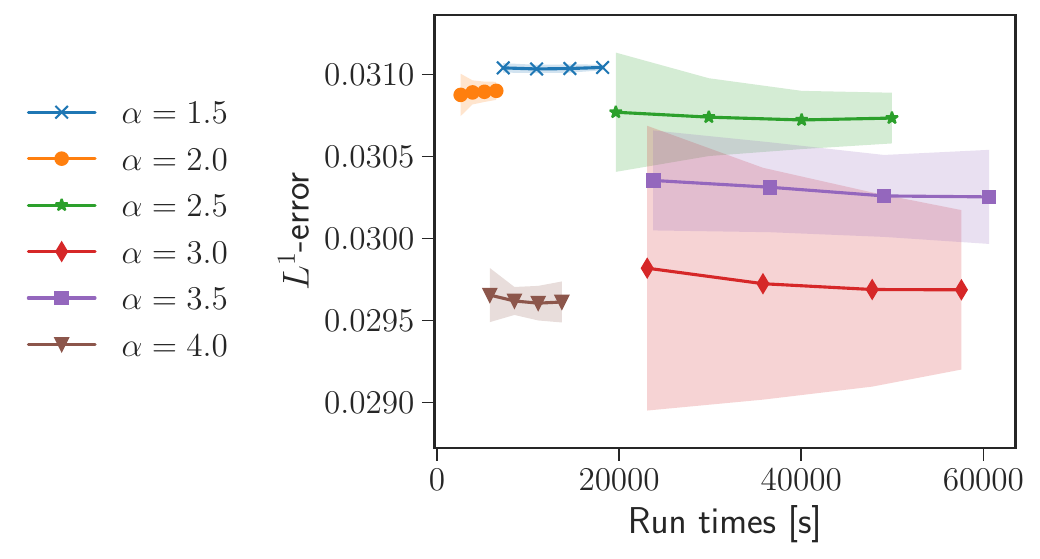}
        \caption*{\hspace{2.25cm}Overview of $L^1$-errors.}
    \end{subfigure}\quad
    \begin{subfigure}[b]{0.4\textwidth}
        \centering
        \includegraphics[height=0.671\linewidth]{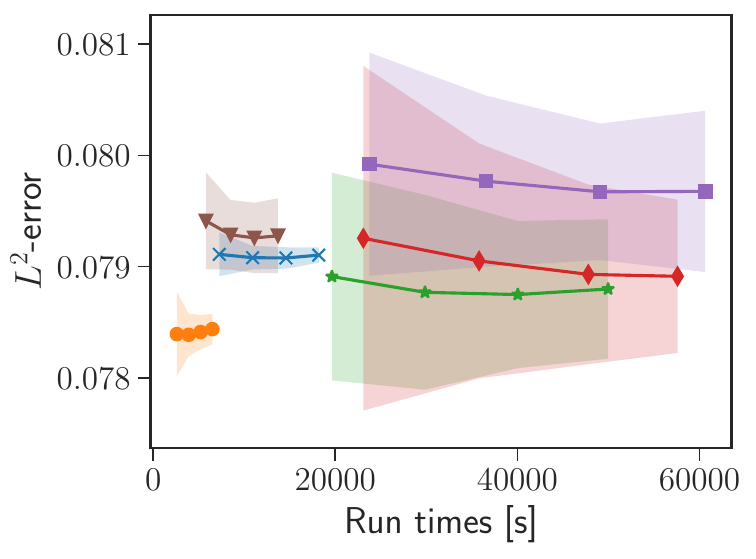}
        \caption*{Overview of $L^2$-errors.}
    \end{subfigure}%
    \\[1em]
    \begin{subfigure}[b]{0.5\textwidth}
        \raggedleft
        \includegraphics[height=0.536625\linewidth]{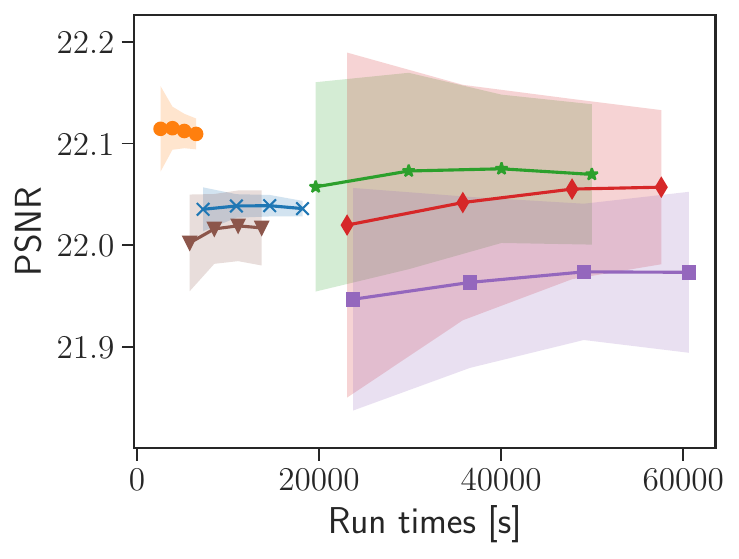}
        \caption*{\hspace{2.25cm}Overview of PSNR values.}
    \end{subfigure}\quad
    \begin{subfigure}[b]{0.4\textwidth}
        \centering
        \includegraphics[height=0.671\linewidth]{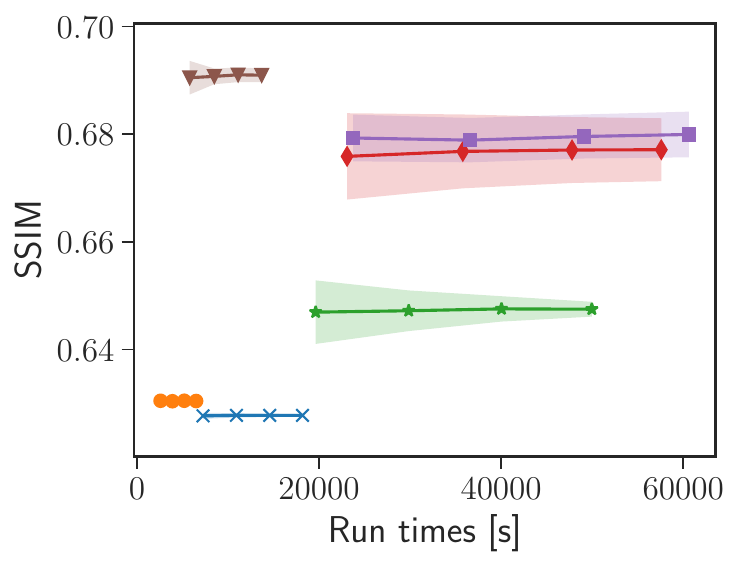}
        \caption*{Overview of SSIM values.}
    \end{subfigure}%
    \caption{Error metrics for ``square/circle'' image and varying $\alpha$ values, as well as iteration counts.}
    \label{fig:error_comparison}
\end{figure}

\paragraph{Cost comparison}
In practice, one is usually not interested in the error of a reconstruction alone, but also in its computational cost.
To be able to provide a thorough comparison, we run the ``square/circle'' experiment for each value of $\alpha$ and to a maximum of $5\times 10^4$ MCMC steps. As before, we discard the first $10^4$ iterations as burn-in. We store the error metrics and associated run times at various MCMC iteration counts during the sampling process. Namely, we record these values at step $2\times 10^4$, $3\times 10^4$, $4\times 10^4$, and $5\times 10^4$ of each individual run. We repeat every experiment a total of 10 times, to compute mean run times as well as mean error metrics.
\review{Additionally, to diagnose convergence in the Markov chains, we compute the effective sample size (ESS), based on the multi-variate potential scale reduction factor (PSRF) introduced in \citep{vats2021revisiting} for multi-variate Markov chains. We note that for every value of $\alpha$, we obtain at least $\sim 5.6 \times 10^4$ effective samples out of $5 \times 10^5$ MCMC samples from combining all 10 chains.
These values indicate that the number of samples we generate is sufficient.
An overview of the obtained ESS values is presented in \Cref{tab:ess} in the appendix.}

The outcomes for the four different error metrics are shown in \Cref{fig:error_comparison}. The shaded areas represent the confidence region of $\pm 2$ standard deviations around the mean, estimated using the 10 experiment repetitions.
A few things stand out in these plots. The cost of each iteration of \Cref{alg:pcn} is considerably smaller than an iteration of \Cref{alg:det_free}. As a result, the experiments that use $\alpha=2$ and $\alpha=4$ have the smallest run times associated with them. The computations for $\alpha=2$ are faster than those for $\alpha=4$ because the precision matrix used has fewer non-zero entries, which leads to sparser Cholesky factors and is exploited by the sparse Cholesky solver that we used. This also explains the large difference in run times between the cases $\alpha=1.5$ and $\alpha=2.5, ~ 3, ~ 3.5$. The preconditioner used in each step of \Cref{alg:det_free} is based on a Cholesky factorisation for the next-largest non-fractional value of $\alpha$. That is, the preconditioner applied in each MCMC step contained less non-zero entries for $\alpha=1.5$ than in the experiments with $\alpha=2.5, ~ 3, ~ 3.5$.

When considering the mean error values obtained for different values of $\alpha$, we observe that $\alpha=2$ seems to perform best w.r.t.\ the $L^2$-error and, equivalently, PSNR. In terms of mean $L^1$-error and SSIM, a value of $\alpha=4$ gives the best result. The value $\alpha=3$, on the other hand, performs well across all four error metrics considered here, and can thus be seen as a good middle ground. However, the confidence region is significantly larger for the values $\alpha=2.5, ~3, ~3.5$, indicating a larger variation in the error metrics over multiple runs.

\begin{figure}
    \centering
    \begin{subfigure}[b]{0.477\textwidth}
        \centering
        \includegraphics[height=0.5625\linewidth]{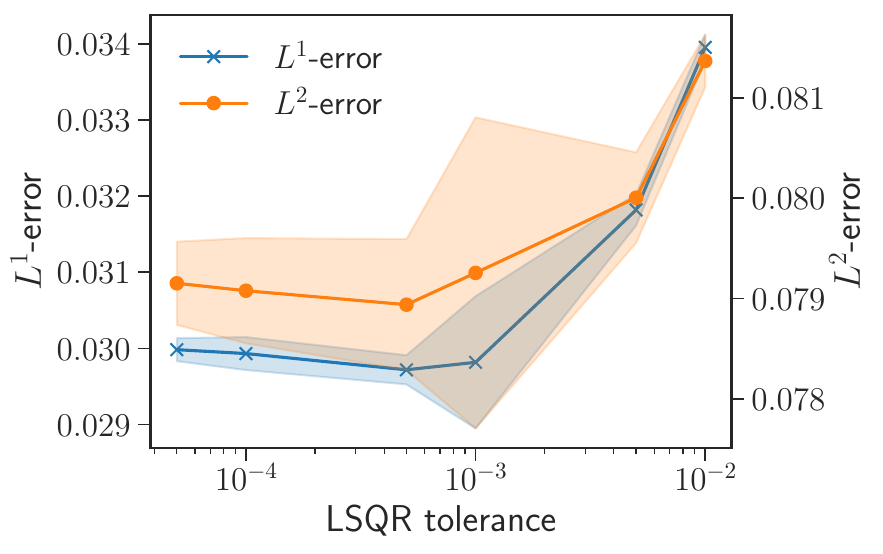}
        \caption*{$L^1$-errors and $L^2$-errors.}
    \end{subfigure}%
    \begin{subfigure}[b]{0.477\textwidth}
        \centering
        \includegraphics[height=0.5625\linewidth]{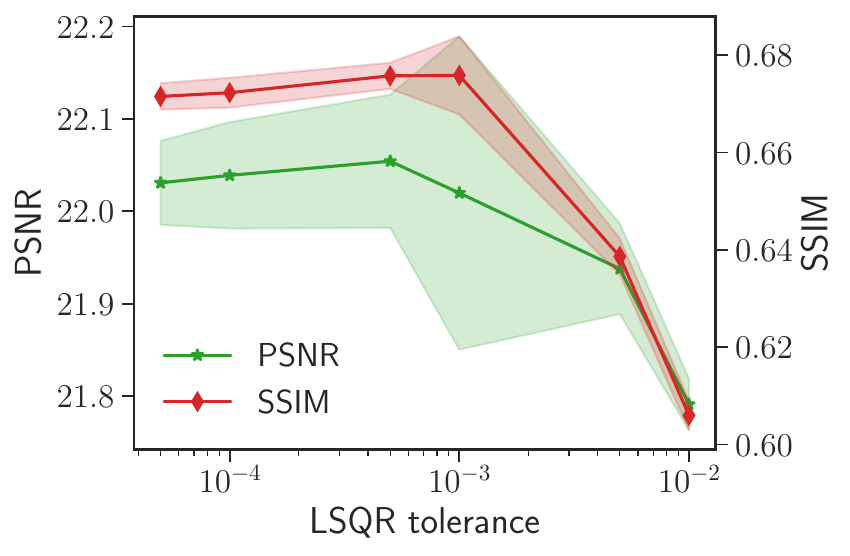}
        \caption*{PSNR and SSIM.}
    \end{subfigure}
    \caption{Error metrics for ``square/circle'' image and varying LSQR tolerances.}
    \label{fig:tol}
\end{figure}

One parameter that greatly influences the computational cost of running \Cref{alg:det_free} is the tolerance set in the LSQR solver. In our implementation of LSQR, we use the stopping rule (S2) for incompatible systems in \citep{paige1982lsqr}, which is based on the system matrix and the residual,. In \Cref{fig:tol}, we compare the errors obtained when reconstructing the ``square/circle'' image for various choices of tolerance parameter. Again, we repeat the experiment a total of 10 times for each tolerance value. The regularity parameter is set to $\alpha=3$. For a large tolerance, we would expect to fail in reconstructing any small-scale features, as terminating the LSQR solver early can be regarded as computing a regularised, smoothed solution \citep{gazzola2020krylov}. For small tolerances, the number of iterations needed to satisfy the stopping criterion increases, leading to a possibly unnecessary increase in computational effort. As we can observe from the plots, it is sufficient to choose a tolerance of about $10^{-3}$ in our example. In our experiments, this tolerance value resulted in roughly 10-15 LSQR iterations per solve. Decreasing this tolerance does not seem to improve the computational outcome.

\subsection{Edge detection}

\begin{figure}
    \centering
    \includegraphics[width=\textwidth]{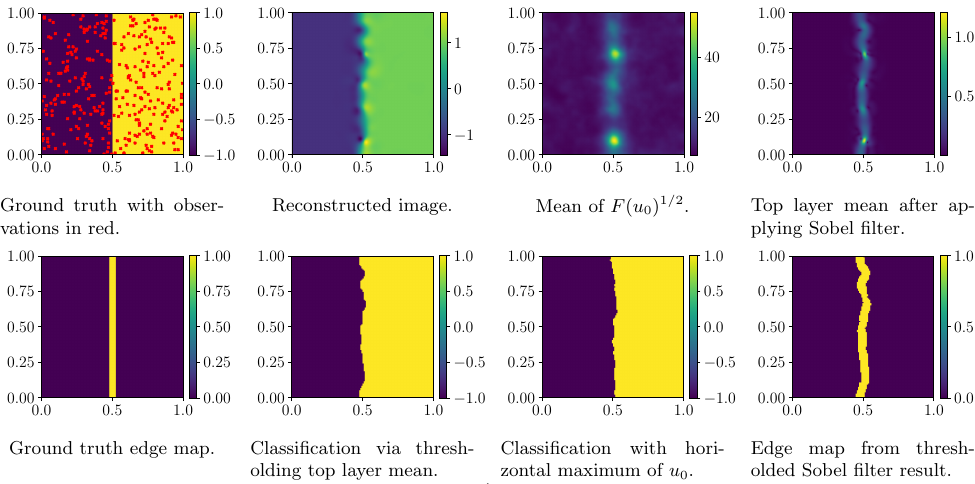}
    \caption{Edge detection results with straight edge and $\alpha=3$.}
    \label{fig:edge_straight}
\end{figure}

In the previous subsection, we were comparing reconstruction results through global error metrics. However, in many applications, local errors are of much larger concern. 
In the following example, we consider the task of finding the location of an edge, given pixel observations at random locations. In such problems, it is not important to represent the image accurately in the entire domain. Instead, only areas around the edge to be found are of consequence, as they determine the outcome of the edge estimation.

Our experiment setup is identical to the one described in \Cref{sec:observation}. However, we now choose the observation locations by uniformly selecting $1/50$ of the $128\times 128$ possible pixels at random. These point observations are not polluted by noise, as the data is binary. We still assume a standard deviation of 0.02 for the noise in our deep GP and stationary GP reconstructions to avoid numerical instability issues.
In what follows, we set the regularity parameter to $\alpha=3$ and use \Cref{alg:det_free} to compute $2\times 10^4$ posterior samples, of which $10^4$ are discarded as burn-in. Another difference in the experiment setup is the ground truth images under consideration. The first such image is shown in the left-most plot in \Cref{fig:edge_straight}: we consider an image divided into a constant left side with value -1 and a constant right side with value 1, so that the edge we are looking for is a straight vertical line in the centre of the image. An image of the ground truth edge (with line thickness 6 pixels) is shown in the bottom left plot in \Cref{fig:edge_straight}. Our approach to finding this edge is to fit a deep GP to the image observations, as in the previous subsection, and to then use the information contained by the deep GP to infer the location of the edge. The result of the reconstruction by deep GP regression is shown in the two upper middle plots in \Cref{fig:edge_straight}.

There are a multiple ways in which we can make use of the fitted deep GP's data. As a first approach, we threshold the estimated mean of the top layer at $0$. That is, every negative value is set to $-1$, and every positive value is set to $1$. The estimated edge is then at the interface of these two regions. We can quantify the quality of this estimation by computing the classification error, that is by calculating the percentage of pixels that have been correctly classified as either negative or positive.
As a second approach, we make use of the information contained in the hidden layer of the deep GP. We note that the average values of $U_0$, or equivalently $F(U_0)$, are large wherever small length scales are advantageous to fit the observations. This is especially the case in the vicinity of edges. In order to exploit this relation, we make use of the simple structure of our ground truth image. In each row, we determine the maximum value of the mean of $F(U_0)$, and then assign the value $-1$ to every pixel to the left of that maximum. The rest is classified as pixels of value $1$. Again, we can assess the classification quality by checking the classification error.
Finally, we estimate the location of the edge through the information contained in the gradient of the top layer mean. That is, we expect large gradient norms in areas close to the edge. In our implementation, we use a Sobel filter \citep{jahne1999principles,scikit-image} to obtain an intensity map of the norm of the top layer gradient. We then threshold this intensity map with a value that is chosen to optimise the F-score \citep{rijsbergen1979information,taha2015metrics} of the resulting edge map. That is, we choose the threshold with the knowledge of the true edge map. In practice, one would usually choose a threshold that gives the best overall performance over a data set and use that value for edge detection in new images. The F-score of our resulting edge map represents how well the edge can be detected with our model.

Examples of all three edge detection approaches are shown in \Cref{fig:edge_straight}. In the bottom row, second plot from the left, we show the classification of the pixels after thresholding the top layer mean. One plot further to the right, we have the classification result relying on finding the maxima of the hidden layer. The bottom right plot shows the edge map obtained from thresholding the gradient intensities obtained through Sobel filtering. The threshold was chosen to optimise the F-score of the edge map with respect to the ground truth edge map in the bottom left plot.

\begin{figure}
    \centering
    \begin{subfigure}[b]{0.477\textwidth}
        \centering
        \includegraphics[height=0.5625\linewidth]{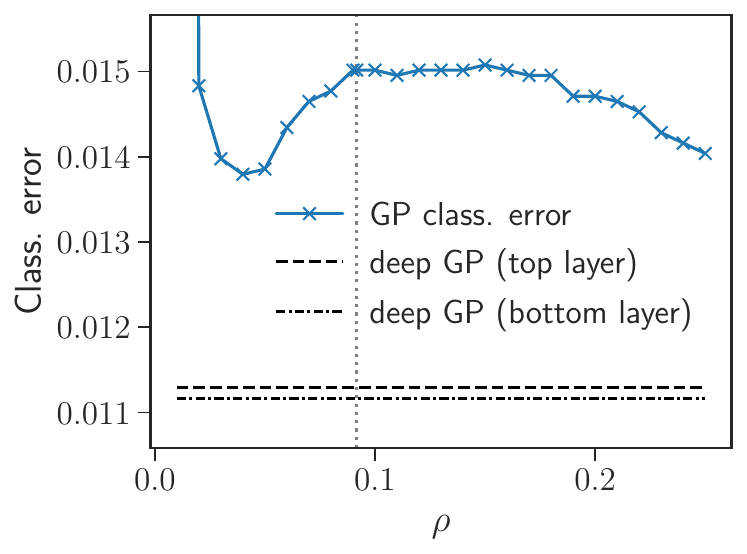}
        \caption*{Classification errors of deep GP and GP results.}
    \end{subfigure}%
    \begin{subfigure}[b]{0.477\textwidth}
        \centering
        \includegraphics[height=0.5625\linewidth]{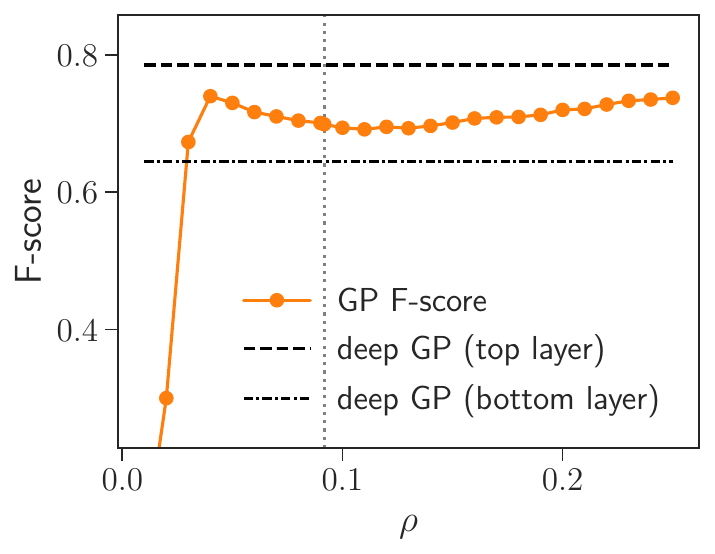}
        \caption*{F-scores of deep GP and GP results.}
    \end{subfigure}
    \caption{Comparison of edge detection errors in experiment with straight edge. }
    \label{fig:edge_straight_errors}
\end{figure}

\begin{figure}
    \centering
    \begin{subfigure}[t]{0.23\textwidth}
        \centering
        \includegraphics[height=0.825\linewidth]{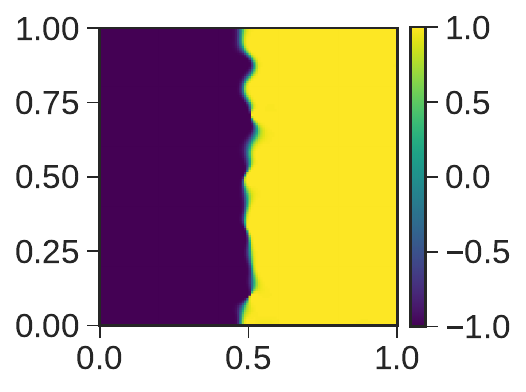}
        \caption*{Mean of top layer samples.}
    \end{subfigure}
    \qquad
    \begin{subfigure}[t]{0.23\textwidth}
        \centering
        \includegraphics[height=0.825\linewidth]{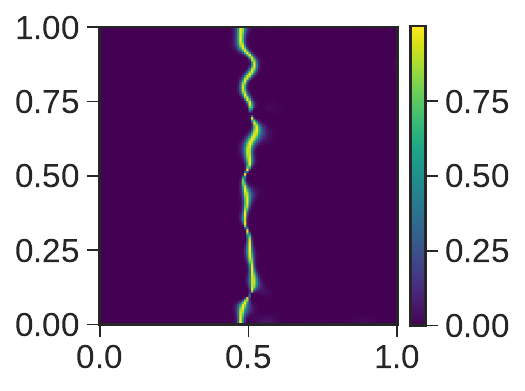}
        \caption*{Point-wise variance values.}
    \end{subfigure}
    \caption{Sample mean and point-wise variances of top layer samples.}
    \label{fig:edge_uq}
\end{figure}

Once more, we compare the results obtained through deep GP regression with what one would achieve with standard, stationary GP regression. We choose the same parameters for our standard GPs as the ones we use in the individual layers of the deep GP. For the length scale, we consider a range of different values. For each of these length scales, we then threshold the stationary GP fitted to the observations and compute the respective classification errors. Analogously to the deep GP top layer, we also apply a Sobel filter to the GP regression results and threshold the resulting intensity maps to obtain F-scores for each length scale value.

We report the results for deep GP and stationary GP edge detection in \Cref{fig:edge_straight_errors}, with the classification errors shown in the left plot, and the F-scores in the right plot. Note that in a perfect reconstruction, we would have a classification error of 0 and an F-score of 1.
The dashed lines indicate the results related to the deep GP-based edge detection making use of information at the top layer, as described earlier. The dash-dotted line in the classification error plot on the left of \Cref{fig:edge_straight_errors} is based on finding the row-wise maxima of the hidden layer. For completeness, the right-hand side plot also includes a dash-dotted line representing the F-score obtained from thresholding the mean of $F(U_0)$, that is, the transformed bottom layer. However, this did not perform well in our experiments. The vertical dotted line represents an estimate of the optimal length scale for standard GP regression. In the reported errors, we can observe that overall, the deep GP performs better than the standard GP models for all tested length scales. This seems to indicate that the improved reconstruction of the ground truth translates into a better detection of the edge as well.
\review{It should be noted that the error metrics we consider here are still global error metrics. The classification error and the F-score are computed over the whole domain, even if we are classifying each pixel in a binary sense. Finding a good (local) metric for edge construction problems is an active research area. Approaches for metrics based on the Sobel filter that aim to measure the accuracy of edge information are given in, e.g., \citep{attar2016image,chen2006edge,martini2012image}.}

\begin{figure}
    \centering
    \includegraphics[width=\textwidth]{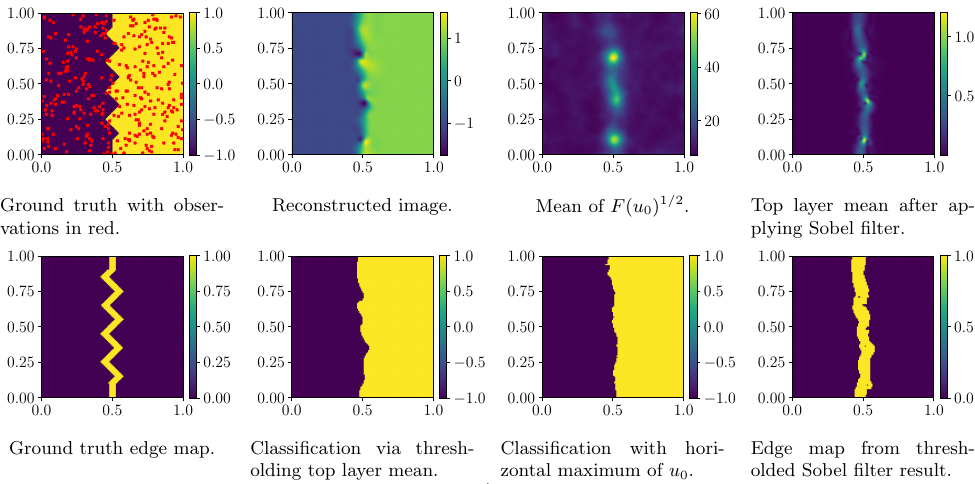}
    \caption{Edge detection results with zigzag edge and $\alpha=3$.}
    \label{fig:edge_zigzag}
\end{figure}

For very large correlation lengths $\rho$, however, we would expect the standard GP to outperform the deep GP. This is related to the way the experiment is set up. As the edge we want to reconstruct is perfectly straight, a large correlation length would help smooth out the reconstructed edge, thus obtaining a better result. However, in applications, we usually do not deal with perfectly straight edges. This motivates the setup of our second edge detection experiment.

\begin{figure}
    \centering
    \begin{subfigure}[b]{0.477\textwidth}
        \centering
        \includegraphics[height=0.5625\linewidth]{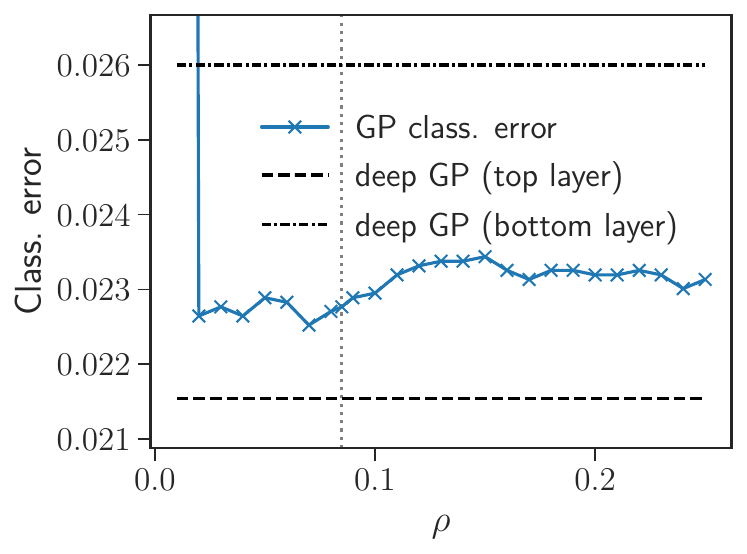}
        \caption*{Classification errors of deep GP and GP results.}
    \end{subfigure}%
    \begin{subfigure}[b]{0.477\textwidth}
        \centering
        \includegraphics[height=0.5625\linewidth]{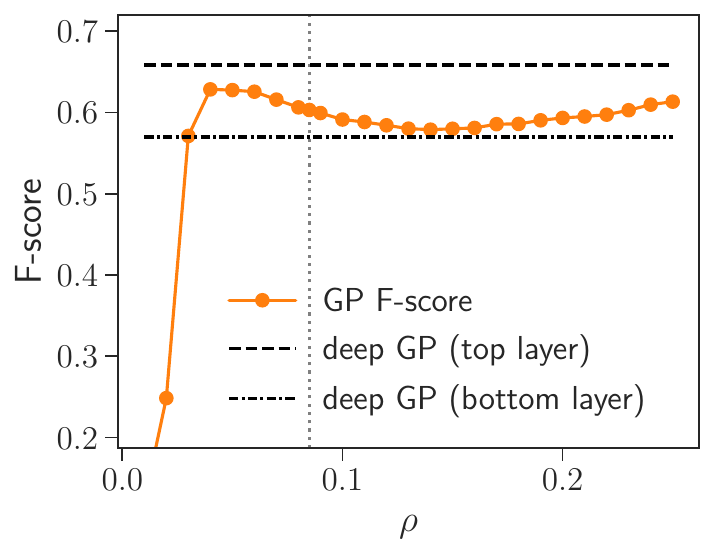}
        \caption*{F-scores of deep GP and GP results.}
    \end{subfigure}
    \caption{Comparison of edge detection errors in experiment with zigzag edge.}
    \label{fig:edge_zigzag_errors}
\end{figure}

We perform the same experiment with a more complex ground truth image, shown in the top left in \Cref{fig:edge_zigzag}. The ground truth edge map is shown in the bottom left.
Most prominently, the edge we want to detect now exhibits small-scale features.
As before, the figure also contains the results of performing deep Gaussian process regression on the observations, as well as the post-processing steps to extract edge information from the fitted deep GP described earlier.
We compare the edge detection results obtained using a deep GP with standard GP regression. The corresponding errors are presented in \Cref{fig:edge_zigzag_errors}.
Overall, the top layer of the deep GP performs the best in both scores. As can be seen in \Cref{fig:edge_zigzag}, the bottom layer of the deep GP does not seem to convey very precise information about the location of the edge. This is also reflected in the error scores shown in \Cref{fig:edge_zigzag_errors}, where the bottom layer is shown to perform badly. For this example, the standard GP achieves similar errors over all length scales under consideration. As smoothing out the reconstructed edge is not advantageous for the ground truth image chosen here, selecting a larger correlation length cannot be expected to improve the obtained error. That is, the edge detection based on a deep GP should outperform the standard GP for any choice of correlation length in this example.

\review{In the edge reconstruction methods outlined above, we did not take into account that we are dealing with stochastic quantities that will have a certain uncertainty associated with them. However, having access to the samples of the hidden layer driving noise allows us to produce samples of both top and hidden deep GP layers. We illustrate this for our edge reconstruction approach relying on thresholding the top layer values. We compute samples of the bottom layer, $u_0$, giving us access to the covariance operator of the top layer $u_1$. Using the formula for the posterior mean and covariance of a Gaussian process \citep{williams2006gaussian}, we then produce a sample of the top layer corresponding to each sample of the bottom layer. As before, we threshold this top layer sample to obtain a classification mapping. The sample mean and point-wise variances of these classification samples are presented in \Cref{fig:edge_uq}. Notably, the variance values indicate that we can be very certain about our classification result when we are far away from the edges. On the reconstructed edge itself, we observe a significant variance. However, we also see that this region of large variance, that is, uncertainty, is very narrow and does not fully cover the true edge we aim at reconstructing. This could mean that some of the deep GP's hyper-parameters, such as the marginal variance, were mis-specified, leading to an over-confident reconstruction.}

\subsection{Image reconstruction with Radon transform}
\begin{figure}
    \centering
    \begin{subfigure}[b]{0.25\textwidth}
        \centering
        \includegraphics[height=0.8\linewidth]{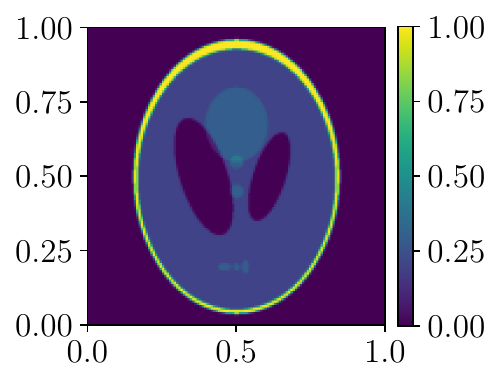}
        \caption*{Shepp-Logan phantom.}
    \end{subfigure}%
    \qquad
    \begin{subfigure}[b]{0.25\textwidth}
        \centering
        \includegraphics[height=0.8\linewidth]{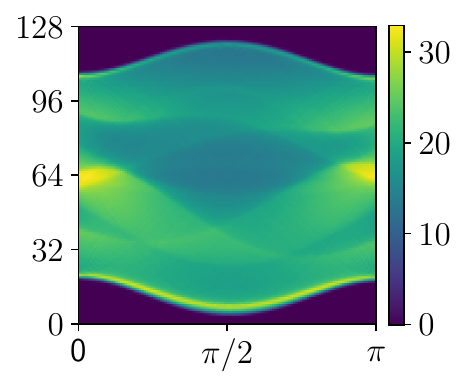}
        \caption*{Full-angle observations.}
    \end{subfigure}%
    \qquad
    \begin{subfigure}[b]{0.25\textwidth}
        \centering
        \includegraphics[height=0.8\linewidth]{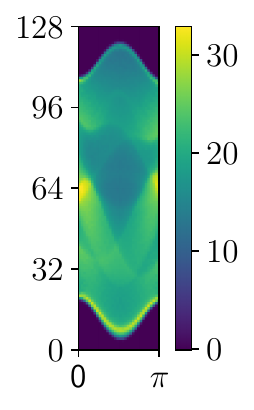}
        \caption*{Sparse-view observations.}
    \end{subfigure}%
    \caption{Ground truth image and observed data for Radon transform experiments.}
    \label{fig:shepp_logan}
\end{figure}

As a final numerical experiment, we consider the inversion of the Radon transform operator \citep{radon1917uber,scikit-image}. This application is a classical problem from medical imaging.
The forward operator in this experiment is the parallel-beam Radon transform as implemented in the scikit-image Python package \citep{scikit-image}. To start with, we consider a full-angle Radon transform with 128 angles chosen uniformly from $[0, \pi)$. As in \Cref{sec:observation}, we assume the presence of Gaussian noise with standard deviation 0.02 polluting the observed data. The deep GP prior is initialised as in the earlier examples, except for the parameter $a$ used in the layer transformation $F$. This is set to $a=400$ in this example. As a ground truth density image, we choose the Shepp-Logan phantom \citep{shepp1974fourier}, rescaled to an image of size $128\times 128$, as shown in the left plot in \Cref{fig:shepp_logan}. The middle plot in the same figure shows the observations obtained by applying the forward operator and adding noise.

\begin{figure}[ht]
    \begin{minipage}{\textwidth}
        \centering
        \begin{subfigure}[b]{0.25\textwidth}
            \centering
            \includegraphics[height=0.8\linewidth]{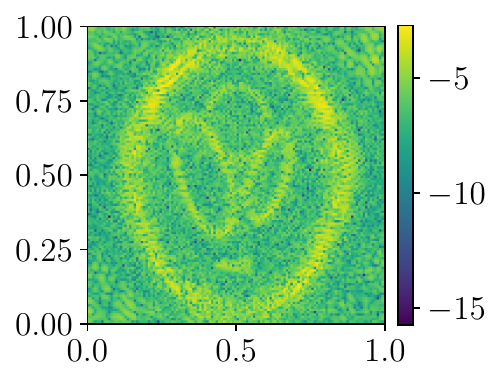}
            \caption*{Point-wise log. errors.}
        \end{subfigure}%
        \qquad
        \begin{subfigure}[b]{0.25\textwidth}
            \centering
            \includegraphics[height=0.8\linewidth]{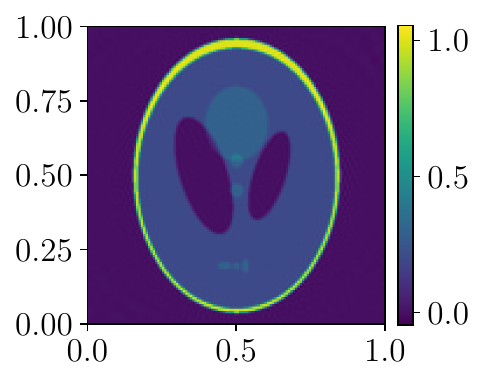}
            \caption*{Reconstructed image.}
        \end{subfigure}%
        \qquad
        \begin{subfigure}[b]{0.25\textwidth}
            \centering
            \includegraphics[height=0.8\linewidth]{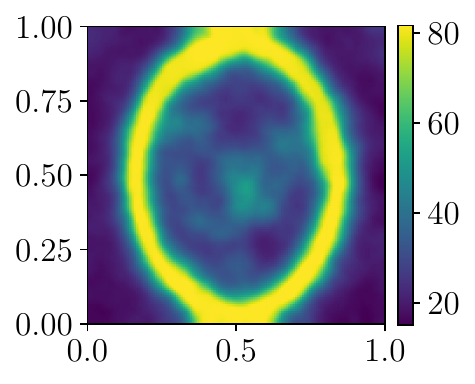}
            \caption*{Estimated mean of $F(u_0)^{1/2}$.}
        \end{subfigure}%
        \caption{\review{Deep GP reconstruction} result for full-angle Radon transform as forward operator and $\alpha=3$.}
        \label{fig:radon}
        \vspace{1em}
        \includegraphics[width=\textwidth]{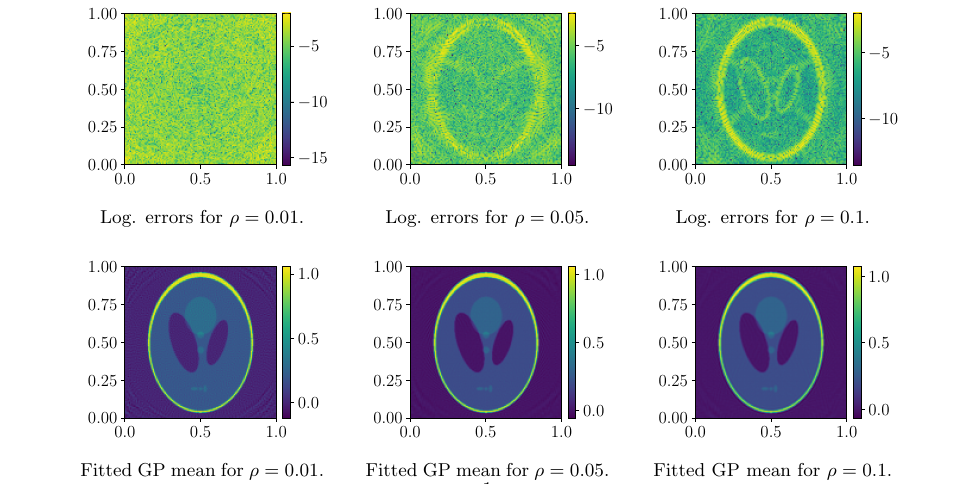}
        \caption{Stationary GP regression results for full-angle Radon transform and $\alpha=3$.}
        \label{fig:radon_gp}
    \end{minipage}
\end{figure}

The Radon transform as a linear operator poses a specific implementational challenge. We mentioned that, in order to use \Cref{alg:pcn}, we employ the Woodbury matrix identity to work with sparse matrices only. While the matrix representation of the Radon transform is sparse, the product $A_h^\ast \Gamma^{-1} A_h$, which appears when evaluating the Woodbury matrix identity, suffers from fill-in. Because this product results in a dense matrix, we cannot use \Cref{alg:pcn} efficiently. Instead, we have to rely on \Cref{alg:det_free} entirely to obtain samples from our posterior distribution. In this experiment, we compute $4\times 10^4$ posterior samples and discard $10^4$ of these as burn-in. In order to speed up the convergence of the LSQR solver used in \Cref{alg:det_free}, we still compute the preconditioner approximating $\Sigma_{L-1}^{-1}$. However, in this example, we only recompute it once every 1000 accepted MCMC steps in order to minimise computation times. The LSQR tolerance parameter is set to $10^{-4}$ for this experiment.

\begin{figure}[p]
    \includegraphics[width=\textwidth]{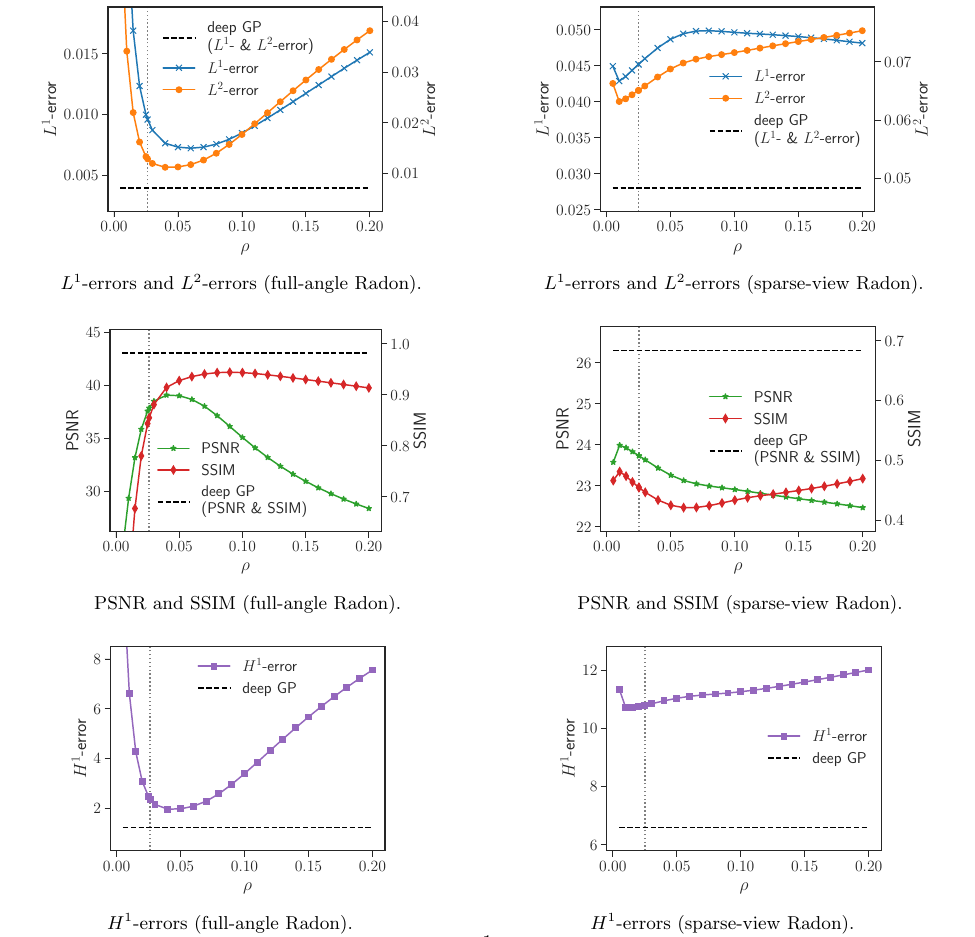}
    \caption{Reconstruction errors for the full-angle Radon transform (left column) and the sparse-view Radon transform (right column) for $\alpha=3$, using standard GP and deep GP priors.}
    \label{fig:radon_errors_a3}
\end{figure}

We present the results for a regularity parameter $\alpha=3$ in \Cref{fig:radon}. As can be observed, the hidden layer picks up the overall structure of the ground truth by assigning a large length scale to the region of the outer circle in the phantom. Some of the smaller scale details on the inside of this circle can be recognised, but are not distinctly visible. Possible reasons for this could be the choice of length scale for the hidden GP layer, or the smoothing effect of the LSQR solver terminating early. Still, using a non-stationary model benefits the reconstruction. This can be seen from the errors obtained for both stationary GP, and deep GP regression, shown in the left column in \Cref{fig:radon_errors_a3}. As earlier, we report the $L^1$- and $L^2$-errors, as well as PSNR and SSIM in two plots. Additionally, we include a plot of the $H^1$-errors, that is, the sum of the $L^2$-error and the $L^2$-norm of the error gradient norm. An estimation of the optimal length scale for standard GP regression is indicated by the vertical dotted line. The plots show that, for all error measures, the reconstruction is improved by using the deep Gaussian process.

We illustrate the influence of the stationary GP's correlation length $\rho$ in \Cref{fig:radon_gp}.
The plots show the GP regression means for $\rho=0.01$, $\rho=0.05$, and $\rho=0.1$ in the bottom row, as well as the associated logarithmic point-wise errors in the top row.
For the smallest length scale, we observe noisy point-wise errors in the entire domain. Still, the edges of the phantom are captured quite well by the reconstruction. This can be explained by the small length scale GP model reconstructing not only small-scale features such as the edges, but also the noise present in the observations. For the larger length scales, the errors caused by the observation noise seem to be less important, and the reconstructed images appear visually smoother.
The error made at the phantom's edges, on the other hand, becomes more and more significant with the correlation length increasing. We can interpret the deep GP results in \Cref{fig:radon} as a combination of the aspects shown in \Cref{fig:radon_gp}: In some areas, the deep GP behaves like a GP with a large length scale, smoothing over the noise that we observe in the reconstructed image. In the areas of the image's main edges, it adopts a small correlation length to recover the edges in a manner similar to the stationary GP with $\rho=0.01$.

\begin{figure}[t]
    \begin{minipage}{\textwidth}
        \centering
        \begin{subfigure}[b]{0.25\textwidth}
            \centering
            \includegraphics[height=0.8\linewidth]{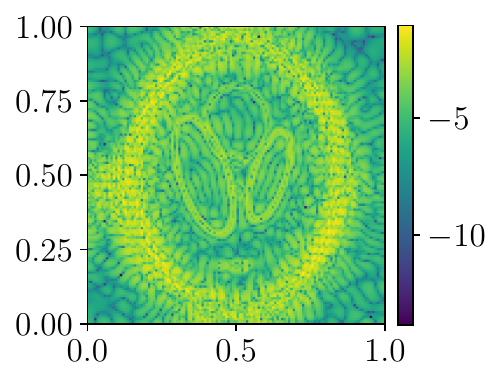}
            \caption*{Point-wise log. errors.}
        \end{subfigure}%
        \qquad
        \begin{subfigure}[b]{0.25\textwidth}
            \centering
            \includegraphics[height=0.8\linewidth]{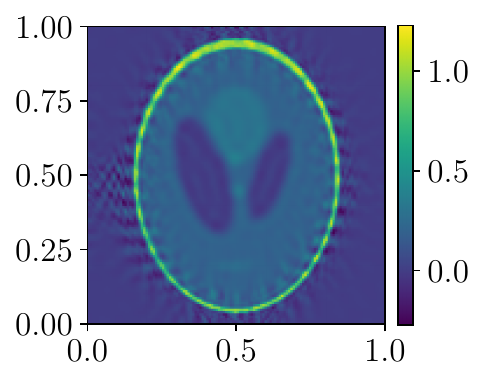}
            \caption*{Reconstructed image.}
        \end{subfigure}%
        \qquad
        \begin{subfigure}[b]{0.25\textwidth}
            \centering
            \includegraphics[height=0.8\linewidth]{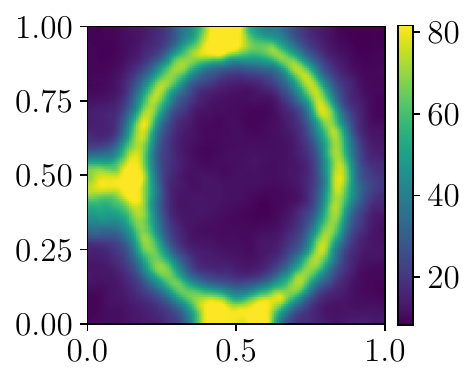}
            \caption*{Estimated mean of $F(u_0)^{1/2}$.}
        \end{subfigure}%
        \caption{\review{Deep GP reconstruction} result for sparse-view Radon transform as forward operator and $\alpha=3$.}
        \label{fig:radon_sparse}
        \vspace{1em}
        \includegraphics[width=\textwidth]{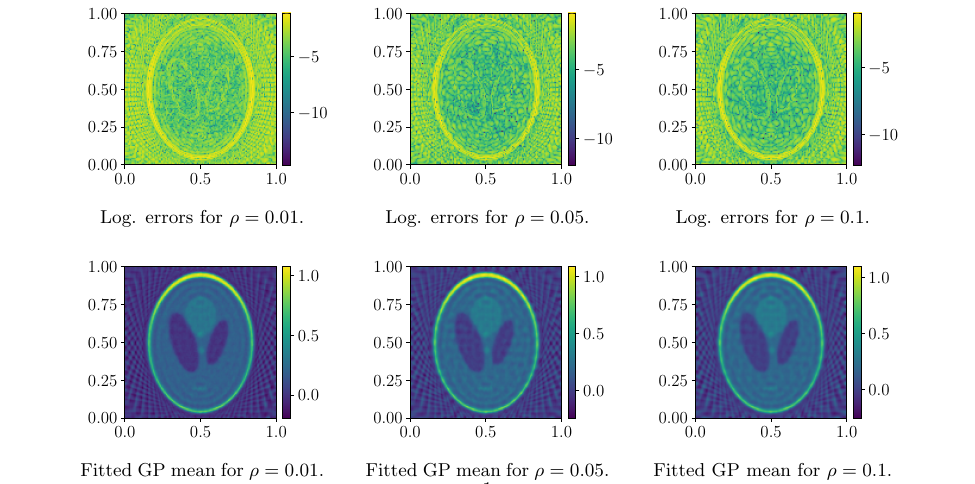}
        \caption{Stationary GP regression results for sparse-view Radon transform and $\alpha=3$.}
        \label{fig:radon_sparse_gp}
    \end{minipage}
\end{figure}

In practice, one usually seeks to reduce the number of X-ray measurements taken. This saves time and keeps the patient's exposure to X-rays to a minimum \citep{frush2003computed,parcero2017impact}. Commonly, this means decreasing the number of angles used in gathering the data. This scenario is mirrored in the following adaptation of our numerical example. As before, we aim at inverting the Radon transform. However, we limit our observations to measurements at 25\%, that is, 32 of the original 128 angles only. These angles are again chosen uniformly from $[0, \pi)$. The observations for this experiment are shown in the right-hand side plot in \Cref{fig:shepp_logan}. Furthermore, we set the parameter $a$ in the layer transformation $F$ to $a=100$, and the LSQR tolerance parameter to $5 \times 10^{-4}$. Otherwise, the setup remains identical to the previous experiment. We present the reconstruction results of this sparse-view example in \Cref{fig:radon_sparse}.
Notably, we are still able to detect the outer edge of the main circular shape from the limited observations in the right-most plot in \Cref{fig:radon_sparse}, depicting the mean of $F(u_0)^{1/2}$.
However, the inner details of the ground truth image are missing now. This is due to a lack of data, as well as the regularising effect of the LSQR solver terminating early. The LSQR tolerance parameter was increased for this sparse-view example as it seemed to improve the convergence of the Markov chain. Furthermore, in the plot of $F(u_0)^{1/2}$, we observe some unnecessarily large values close to the edges. These artifacts are most likely caused by our choice of homogeneous Neumann boundary conditions. The effect of the artifacts is clearly visible as areas of increased noise in the plot on the left in \Cref{fig:radon_sparse}, showing the point-wise reconstruction errors.

Once more, we compare the deep GP reconstruction with stationary GP regression. The results for three selected length scales are shown in \Cref{fig:radon_sparse_gp}. The effect of the length scale on the reconstruction results is similar to the full-view observation case in \Cref{fig:radon_gp}. Overall, however, the noise is much more prominently visible over all three length scales considered. Once more, the deep GP reconstruction has the advantage of employing locally varying length scales. This is also reflected in the corresponding error metrics, shown in the right column plots in \Cref{fig:radon_errors_a3}. Here, we observe a clear margin of improvement over the stationary GP regression results for all error metrics.

\section{Conclusion}
\label{sec:conclusion}
In the preceding sections, we have described a computational framework for performing image reconstruction with deep Gaussian process priors. Notably, we do not restrict the value of the smoothness parameter $\alpha$ to be an even number, but rather allow any admissible value $\alpha > d/2$.
\review{Being able to choose the smoothness parameter freely can be an advantage when the true regularity of the ground truth is known, e.g. in applications such as \citep{minasny2005matern,sang2011covariance}. Furthermore, our approach gives the modeller the flexibility to test and try both fractional and non-fractional values of the smoothness parameter when its true value is not known.}
The methodology used in the above combines a number of computational tools from different areas of computational mathematics and statistics.
Firstly, we make use of the stochastic PDE representation for \matern-type Gaussian processes. For values of $\alpha$ that lead to a fractional SPDE, we employ rational approximation to represent the covariance of the corresponding Gaussian processes. As we employ our deep GP model as the prior distribution in a Bayesian inverse problem, we furthermore use MCMC-type algorithms to obtain samples of the resulting posterior. Depending on the value of $\alpha$ and the sparsity of the forward operator, we propose to either use a pCN-MCMC algorithm when possible, or a determinant-free MCMC algorithm with pCN proposals. In the second case, we also briefly described a preconditioner to speed up the iterative LSQR solver employed in the potential evaluation.

We then tested our approach in a series of numerical experiments. In these experiments, we solve image reconstruction tasks using a deep Gaussian process prior. Throughout, we could observe that the use of a non-stationary model was beneficial in reconstructing the ground truth signal, as compared to a stationary Gaussian process prior. Furthermore, our computational framework allowed us to directly compare between results for fractional and non-fractional values of the smoothness parameter.

In future work, we plan to explore further possibilities to evaluate the posterior. One such approach could be to optimise the posterior density directly to obtain a maximum a posteriori estimate. This could lead to a faster solution algorithm, while losing the ability to explore the posterior through sampling.

\review{Since we are considering numerical solutions to (stochastic) PDEs with locally varying length scales in combination with the finite element method, it appears natural that one could adapt the finite element mesh to the current local length scales. That is, one could refine the mesh in areas where the local length scale is small to better capture small-scale features of the solution, and coarsen the mesh in areas of large local length scales, reducing computational cost. We did not explore this idea further in this work, since re-meshing and re-computing the associated stiffness and mass matrices introduces significant additional computational cost. However, there might be a clever way to mitigate these costs, e.g., by ensuring the mesh is adapted only when one has detected that this has sufficient computational benefit in terms of computation speed or accuracy.}

A different route to be investigated is to expand the deep GP model to incorporate directional information \citep{pragliola2023total,llamazares2024parameterization}. The deep GP we use is isotropic by construction, whereas spatial data is often intrinsically anisotropic. One application could be climate modelling in a region containing mountain ranges.
So far, we have also restricted ourselves to \matern-type GPs. The reason for this is the stochastic PDE representation, allowing for the efficient approximation of the covariance via sparse matrices. To make use of deep Gaussian processes with a different covariance structure would require new approaches for sampling the posterior distribution. This could be beneficial to applications that are not well captured by \matern-type Gaussian processes.

\section*{Acknowledgements}
The authors acknowledge the use of the Heriot-Watt University high-performance computing facility (DMOG) and associated support services in the completion of this work. The authors would like to thank the Isaac Newton Institute for Mathematical Sciences for support and hospitality during the programme ``The mathematical and statistical foundation of future data-driven engineering'' when work on this paper was undertaken (EPSRC grant EP/R014604/1). ALT was partially supported by EPSRC grants EP/X01259X/1 and EP/Y028783/1.

\appendix
\section{Additional results}
\label{sec:add_results}

In this section, we include results for experiments with an additional ground truth image, shown on the right in \Cref{fig:ground_truth}.
In \Cref{sec:observation}, we presented two examples where the ability of the deep GP prior to represent multi-scale features improved the overall reconstruction result. Specifically, we observed better errors when comparing to stationary GP regression over a range of correlation length values.
However, using a deep GP prior does not automatically lead to a more accurate reconstruction. We present an example here to illustrate the point.

\begin{figure}[ht]
    \centering
    \includegraphics[width=\textwidth]{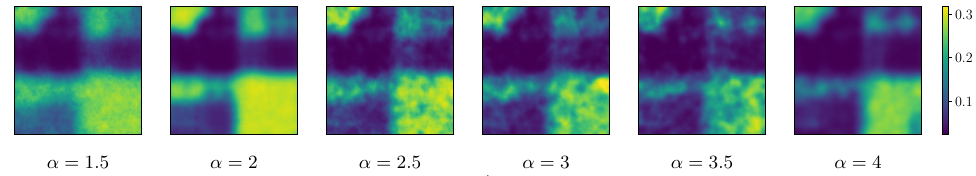}
    \caption{Correlation lengths used in deep GP reconstructions of ``multiscale'' image with varying $\alpha$.}
    \label{fig:corr_lengths_multiscale}
\end{figure}

The setup is the same as described earlier in \Cref{sec:observation}. We show the length scales estimated using \Cref{alg:pcn,alg:det_free} in \Cref{fig:corr_lengths_multiscale}. The overall structure of the ground truth image is clearly visible. However, the image is too complex for all the details to show in the estimation of the bottom GP layer. Especially the small square shapes are not represented well. This is likely due both to the fixed correlation length of the bottom layer, as well as the low spatial density of the observation data points.

\begin{figure}[ht]
    \centering
    \begin{subfigure}[b]{0.413\textwidth}
        \centering
        \includegraphics[height=0.65\linewidth]{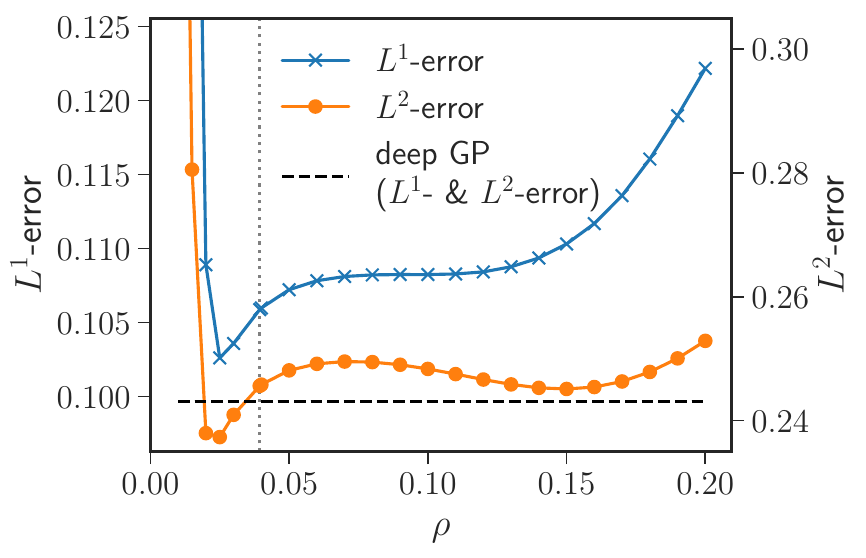}
        \caption*{$L^1$-errors and $L^2$-errors.}
    \end{subfigure}%
    \qquad
    \begin{subfigure}[b]{0.413\textwidth}
        \centering
        \includegraphics[height=0.65\linewidth]{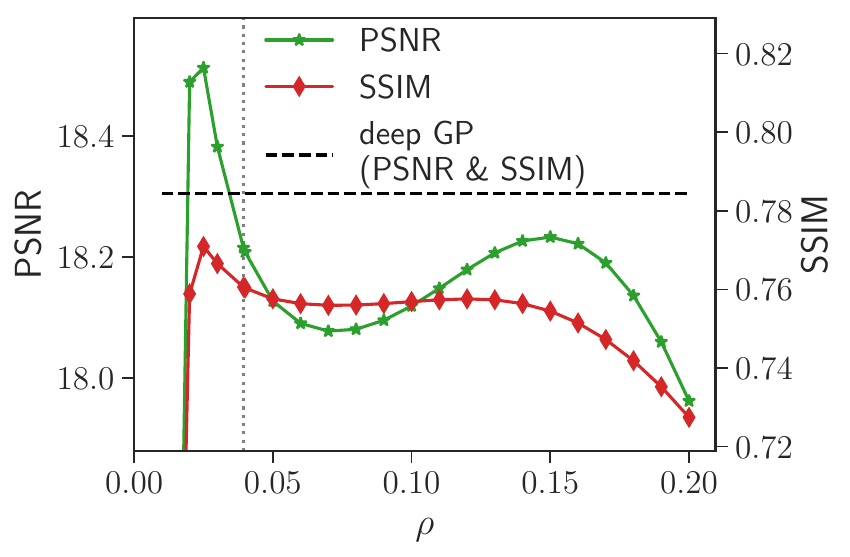}
        \caption*{PSNR and SSIM.}
    \end{subfigure}%
    \caption{``Multiscale'' ground truth image and reconstruction errors for $\alpha=3$, using standard GP and deep GP priors.}
    \label{fig:multiscale_errors_a3}
\end{figure}

For the specific choice of $\alpha=3$, we once more compare the performance of the deep GP reconstruction with standard GP regression over a range of different correlation lengths. The results are shown in \Cref{fig:multiscale_errors_a3}.
While for most stationary correlation length values, including the optimal length scale found by maximum likelihood estimation, the deep GP prior gives better reconstruction errors, there are length scales for which standard GP regression performs better in terms of $L^2$-error/PSNR. We see this as an example where the observational data is not enough to recover a meaningful length scale estimate. As a result, the deep GP prior's flexibility in representing multiple length scales does not have the expected benefit on the reconstruction of the image.

\review{Furthermore, we include a table, \Cref{tab:ess}, showing the effective sample size (ESS) estimates we obtained for the experiment comparing run times for different values of $\alpha$ in \Cref{sec:observation}. The first row of values gives the average over 10 chains, where the ESS is evaluated for each chain individually. The second row gives the ESS of all 10 Markov chains combined.}

\begin{table}[ht]
    \centering
    \begin{tabular}{@{}lrrrrrr@{}} \toprule
        Regularity parameter $\alpha$ & 1.5 & 2 & 2.5 & 3 & 3.5 & 4 \\ \midrule
        Avg. single-chain ESS & 1\,749 & 1\,126 & 252 & 183 & 168 & 635 \\
        Multi-chain ESS       & 156\,483 & 126\,163 & 67\,195 & 57\,390 & 56\,423 & 96\,326 \\ \bottomrule
    \end{tabular}
    \caption{Effective sample sizes for cost comparison experiment in \Cref{sec:observation}, rounded to the closest integer.}
    \label{tab:ess}
\end{table}

\bibliographystyle{siam}
\bibliography{references}

\end{document}